\documentclass[11pt,a4paper]{article}
\usepackage[utf8]{inputenc}
\usepackage[english]{babel}
\usepackage{amsmath}
\usepackage{amsfonts}
\usepackage{amssymb}
\usepackage{graphicx}
\usepackage{algorithmic}
\usepackage{algorithm2e}
\usepackage{a4wide}
\usepackage{verbatim}
\RequirePackage{ifthen}
\RequirePackage[T1]{fontenc}
\RequirePackage{amsthm}
\RequirePackage{color}
\RequirePackage{listings} 
\RequirePackage{epstopdf}
\RequirePackage{emptypage}
\RequirePackage[nolist]{acronym}
\RequirePackage{multirow}
\RequirePackage{enumerate,enumitem}
\RequirePackage{caption}
\RequirePackage{tikz}
\usetikzlibrary{arrows.meta}
\RequirePackage{latexsym,amscd,makeidx}
\RequirePackage{exscale}
\RequirePackage{mathrsfs}
\RequirePackage{hyperref}
\RequirePackage{url}
\RequirePackage{booktabs}
\RequirePackage{esdiff}
\RequirePackage{listings}
\RequirePackage{tabularx}
\newcolumntype{x}[1]{!{\centering\arraybackslash\vrule width #1}}
\RequirePackage[numbers]{natbib}
\RequirePackage{verbatim}
\RequirePackage{comment}
\usepackage{mathtools}
\RequirePackage{cleveref}
\crefname{equation}{}{}
\Crefname{Equation}{}{}
\usepackage{pgfplots} 
\usepackage{pgfgantt}
\usepackage{pdflscape}
\usepackage[pagewise]{lineno}
\pgfplotsset{compat=newest} 
\pgfplotsset{plot coordinates/math parser=false}
\numberwithin{equation}{section}
	\newif\ifshowcomments								
	\showcommentsfalse				

	\definecolor{mygrey}{gray}{0.65}
	
	\ifshowcomments
		
	\else
	 	\excludecomment{bcomment}
	\fi
	
	\ifshowcomments
		\newcommand{\scomment}[1]{{\color{mygrey} [#1]}}  
	\else
		\newcommand{\scomment}[1]{\iffalse #1 \fi}
	\fi

\newlength\fwidth

\theoremstyle{plain}						
\newtheorem{theorem}{Theorem}[section]

\theoremstyle{definition}
\newtheorem{definition}[theorem]{Definition}

\newtheorem{assumption}[theorem]{Assumption}
\theoremstyle{remark}
\newtheorem{remark}[theorem]{Remark}

\newcommand{\N}{\mathbb{N}}
\newcommand{\R}{\mathbb{R}}

\newcommand{\fpustar}{f^{-1}_+}

\newcommand{\rv}[1]{{#1}}

\graphicspath{{Graphik/}}
\begin{document}

\title{Conservation laws with discontinuous flux function on networks: a splitting algorithm}
\author{Jan Friedrich\footnotemark[1], \; Simone Göttlich\footnotemark[2], \; Annika Uphoff\footnotemark[2]}

\date{\today}

\maketitle
\footnotetext[1]{RWTH Aachen University, Institute of Applied Mathematics, 52064 Aachen, Germany (friedrich@igpm.rwth-aachen.de).}
\footnotetext[2]{University of Mannheim, Department of Mathematics, 68131 Mannheim, Germany (goettlich@uni-mannheim.de).}

\begin{abstract}
In this article, we present an extension of the splitting algorithm proposed in \cite{Towers} to networks of conservation laws with \rv{piecewise linear} discontinuous flux functions in the unknown. We start with the discussion of a suitable Riemann solver at the junction and then describe a strategy how to use the splitting algorithm on the network. \rv{In particular, we focus on two types of junctions, i.e., junctions where the number of outgoing roads does not exceed the number of incoming roads (dispersing type) and junctions with two incoming and one outgoing road (merging type).} Finally, numerical examples demonstrate the accuracy of the splitting algorithm by comparisons  to the exact solution and other approaches used in the literature.

  \medskip


  \medskip

\end{abstract}

\section{Introduction}
Traffic flow models based on scalar conservation laws with continuous flux functions are widely used in the literature. For a general presentation of the models, we refer to the books \cite{garavellohanpiccoli2016book,garavello2006traffic,Treiber2013} and the references therein. 
Extensions to road traffic networks have been also established. We mention in particular the contributions  \cite{coclite2005,holden1995}, where the authors introduce the coupled network problem and show the existence of solutions.   
Within this article, we are concerned with the special case of scalar conservation laws with discontinuous flux in the unknown that are motivated in the traffic flow theory by the observation of a gap between the free flow and the congested flow regime \cite{Ceder1976ADT,ceder1976further,Edie1961}. This phenomenon generates an interesting dynamical  behavior called {\em zero waves}, i.e. waves with infinite (negative) speed but zero wave strength, and has been investigated in recent years either from a theoretical or numerical point of view, see for instance \cite{BurgerChalonsVilladaMulticlass, Lu2009,ImplicitMethod,Towers,Wiensa} or more generally \cite{bulicek2011,Carillo2003,Dias2005,gimse1993}.

To the best of our knowledge, the study of scalar conservation laws with discontinuous flux functions on networks is still missing in the traffic flow literature. However, in the context of supply chains with discontinuous flux such considerations have been already done \cite{festa2019,gottlich2013discontinuous}. We remark that supply chain models differ essentially from traffic flow models due to simpler dynamics and different coupling conditions.

In this work, we aim to derive a traffic network model, where the dynamics on each road are governed by a scalar conservation law with discontinuous flux function in the unknown. \rv{For simplicity, we restrict to piecewise linear flux functions.} Special emphasis is put on the coupling at junction points to ensure a unique admissible weak solution. \rv{In particular, we focus on dispersing junctions where the number of incoming roads does not exceed the number of outgoing roads and merging with two incoming and one outgoing road. The latter type of junction can be extended to the case of multiple incoming roads and a single outgoing one.} In order to construct a suitable numerical scheme that is not based on regularization techniques we adapt the splitting algorithm originally introduced in \cite{Towers}. Therein, the discontinuous flux is decomposed into a Lipschitz continuous flux and a Heaviside flux such that a two-point monotone flux scheme, e.g Godunov, can be employed in an appropriate manner. \rv{This algorithm has been studied in \cite{Towers} for the case of a single road only. However, in the network case, multiple roads with possibly disjunctive flux functions need to be considered at a junction point to ensure mass conservation. Hence, the key challenge is to determine the correct flux through the junction in an appropriate manner}.  Therefore, a detailed case distinction in accordance with the theoretical investigations is provided for the different types of junctions. The numerical results validate the proposed algorithm for some relevant network problems.

The paper is organized as follows: in section \ref{sec:disc} we discuss the basic model and Riemann problems which permit to derive an exact solution. 
We extend the modeling framework to networks in section \ref{sec:discnetworks} and focus on the coupling conditions. In section \ref{sec:splitnetworks}, we introduce how the splitting algorithm \cite{Towers} can be extended to also deal with the different types of junctions. Finally, we present a suitable discretization and numerical simulations in section \ref{sec:num}.

\section{LWR model with discontinuous flux on a single edge}
\label{sec:disc}
In this section, we briefly recall the case of the Lighthill-Whitham-Richards (LWR) model \cite{Lighthill1955,Richards1956} on a single edge with a flux function having a single decreasing jump at $ u^* \in [0, u_{\text{max}}]$ that has been intensively studied in e.g. \cite{Lu2009,Towers,Wiensa}.
The shape of this flux function is inspired by empirical research which suggest that the fundamental diagram is of a reverse $\lambda$ shape and consist of a free flow area and a drop to a congested area, see \cite{Ceder1976ADT,ceder1976further,Edie1961}. 

\begin{figure}
  \centering
  \includegraphics[width=0.5\textwidth]{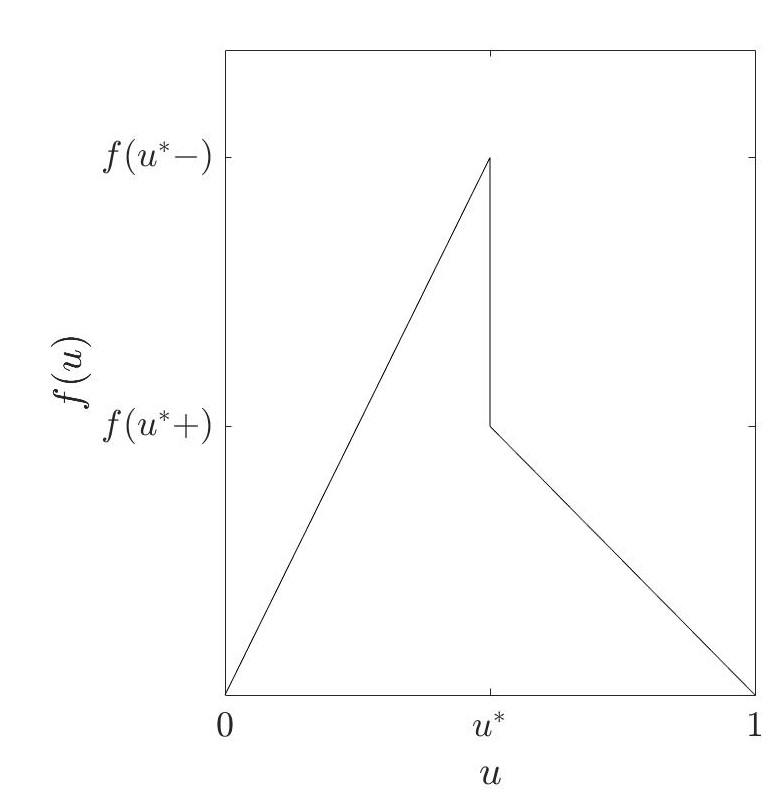}
   \caption[Example for discontinuous flux functions.]{\rv{An example for a piecewise linear discontinuous flux function with discontinuity at $u^*$, cf. \cite{Lu2009}.}}
   \label{fig:DiscFlux}
\end{figure}

Following \cite{Towers}, we consider the scalar conservation law
\begin{equation}\label{ScalarConservationLaw}
\begin{cases}
u_t + f(u)_x = 0, & (x,t) \in (a,b) \times (0,T) =: \Pi_T, \\
u(x,0) = u_0(x) \in [0, u_{\text{max}}], & x \in (a,b),  \\
u(a,t) = r(t) \in [0, u_{\text{max}}], & t \in (0,T), \\
u(b,t) = s(t) \in [0, u_{\text{max}}], & \mathcal{F}(t) \in \tilde f (s(t)), t \in (0,T).
\end{cases}
\end{equation}
More precisely the flux function is defined as follows, 
\begin{equation}\label{FluxFunction}
f(u) = \begin{cases} f_1(u) \ \hat= \text{ Flux 1},& \text{ if } u \in [0,u^*],\\ 
f_2(u) \ \hat= \text{ Flux 2},& \text{ if }u \in (u^*,u_\text{max}]. \end{cases}
\end{equation}
We denote $f_1(u^*) = f(u^*-)$ and $f_2(u^*) = f(u^*+)$, where $f(u^*-) > f(u^*+)$.
In order to assign a single numerical value to $f(u^*)$ we arbitrarily define $f(u^*) = f(u^*-)$.
\rv{One example is} shown in figure \ref{fig:DiscFlux}.
The jump has magnitude 
\begin{equation}
\alpha := f(u^* - ) - f(u^* +).
\end{equation}
As usual we require for the flux function $f : [0, u_\text{max}] \rightarrow \R^+$ with $f(0)=f(u_{\text{max}}) = 0$ and in addition,
\rv{ we assume $f_1$ to be linear increasing and $f_2$ linear decreasing.}

As in \cite{Towers}, the multivalued version of $f$ is defined by
\begin{align}
\tilde f(u) = \begin{cases}
f(u), & u \in [0, u^*),\\
[f(u^*+), f(u^*-)], & u = u^*, \\
f(u), & u \in (u^*, u_{\text{max}}].
\end{cases}
\end{align}

Finally, we have to discuss the imposed boundary conditions at $x=b$ in \eqref{ScalarConservationLaw}.
Due to the flux discontinuity we have a non-standard boundary condition $\mathcal{F}(t) \in \tilde f (s(t))$.
If the boundary condition at the right is $s(t) \neq u^*$, then the additional boundary condition is redundant, i.e., $\mathcal{F}(t)=f(s(t))$.
Otherwise, if $s(t)=u^*$, we have to choose between two possible flux values $\tilde f(u^*) \in \{f(u^*+), f(u^*-)\}$. 
Then, the function $\mathcal{F}$ assigns a flux value to $u(b,t) = s(t)$ according to \cite{Towers} 
\begin{equation}\label{eq:fluxboundary}
\mathcal{F}(t) = \begin{cases}
f(u^* - ), & \text{if the traffic ahead of } x = b \text{ is free-flowing,}\\
f(u^* + ), & \text{if the traffic ahead of } x = b \text{ is congested}.
\end{cases}
\end{equation}
The state of traffic ahead of $ x=b$ is additional information that is determined independently from the other data of the problem. 
\begin{remark}\cite[Remark 1.3]{Towers}
We note that for the boundary condition at the left end the state of traffic ahead of $ x = a$ is known from the start values.
Hence, no additional boundary conditions are necessary.
\end{remark}
The following assumptions are important for the proof of existence and uniqueness of solutions.
\begin{assumption}\label{assumptionSplit}\cite[Assumption 1.1]{Towers}
The initial data satisfies $u_0(x) \in [0, u_{\text{max}}]$ for $x \in (a,b)$ and $u_0 \in \text{BV}([a,b]).$ The boundary data satisfies $r(t),s(t) \in [0, u_{\text{max}}]$ for $ t \in [0,T]$, and $r,s \in \text{BV}([0,T])$.
\end{assumption}
A weak solution is intended in the following sense:
\begin{definition}\cite[Definition 1.1]{Towers}\label{WeakSolutionDisc}
A function $u \in L^{\infty} (\Pi_T)$ is called a weak solution to the initial-boundary value problem \eqref{ScalarConservationLaw} if there exists a function $v \in L^{\infty}(\Pi_T)$ satisfying $$ v(x,t) \in \tilde f(u(x,t)) \text{ a.e.}$$ such that for each $\psi \in C_0^1((a,b) \times[0,T)), $
$$ \int_0^T \int_a^b (u\psi_t + v\psi_x)\text{d}x\text{d}t + \int_a^bu_0(x)\psi(x,0)\text{d}x = 0.$$
\end{definition}

As usual, weak solutions do not lead to a unique solution and additional criteria are necessary to rule out physically incorrect solutions.
In particular, the discontinuity of the flux prohibits from directly using the classical approaches.
Note that in \cite{Towers} an adapted version of Oleinik's entropy condition \cite{Evans} is used to single out the correct solution, while in \cite{Wiensa} the convex hull construction \cite[Chapter 16]{Leveque2002} is used to construct solutions to Riemann problems. 

Here, we will concentrate on the convex hull construction.
For completeness we will shortly recall the solutions to Riemann problems considered in \cite{Wiensa} as they are essential in order to construct a Riemann solver at a junction.

\subsection{Riemann solutions}\label{sec_ZeroWaveRP}

We consider a Riemann problem with initial data $(u_L,u_R)$ and its exact solution.
For simplicity we will focus on a piecewise linear flux function as in \cite{Wiensa} and in figure \ref{fig:DiscFlux}.
The exact solution is then derived with a regularization of the flux and the convex hull construction.
\rv{As regularizations are not unique we choose a simpler regularization than the one in \cite{Wiensa}. 
In fact, we use the regularization introduced in \cite{Lu2009} for a piecewise quadratic flux function and apply it to the piecewise linear flux function.
In the following we simply recall the results from \cite{Lu2009} in case of a linear flux function for the sake of completeness and in order to understand better the upcoming analysis.
Note that the the results also coincide with the ones of \cite{Wiensa}.}

\begin{figure}[htb]
  \centering
  \includegraphics[width=0.7\textwidth]{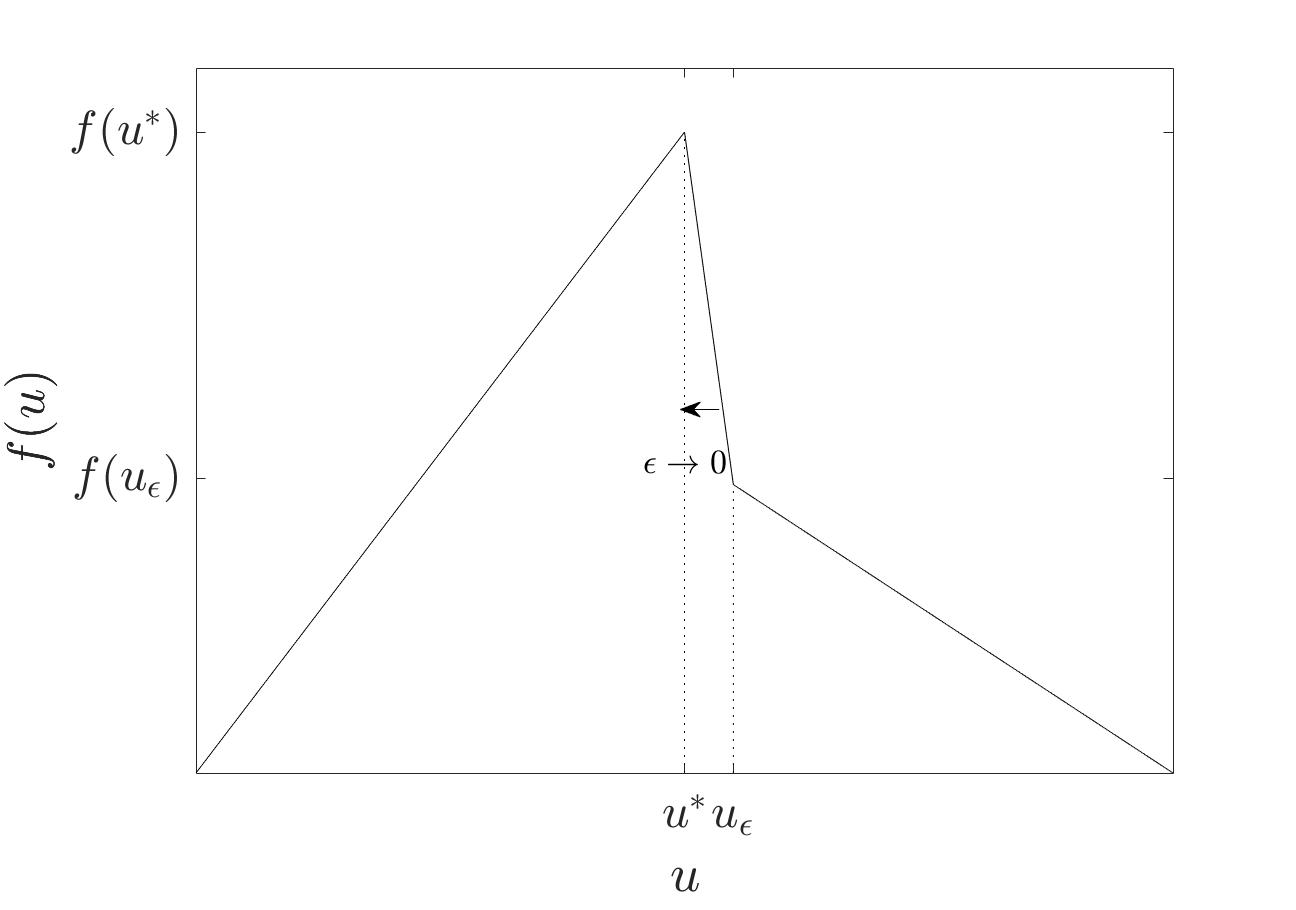}
   \caption[Regularized discontinuous flux function.]{ Example of the regularization of a piecewise linear flux function.}
   \label{fig:Linearization}
\end{figure}

\begin{figure}[htb]
  \centering
  \includegraphics[width=0.5\textwidth]{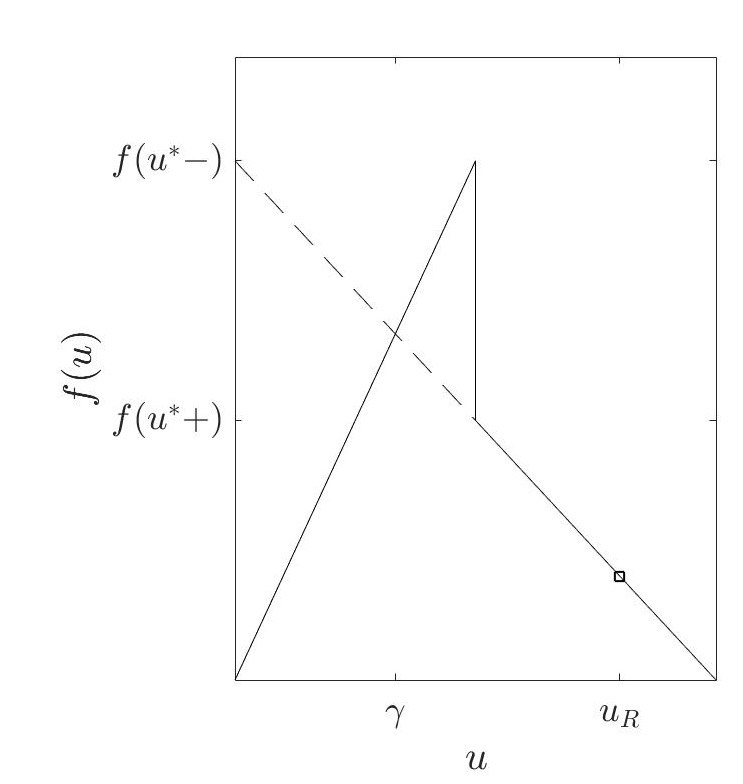}
   \caption[Convex Hull construction - function $h$]{ The dotted line shows the function $h$, connecting $(u^*,f(u^*+)) $ and $u_R$. The value $\gamma$ is the intersection point. }
   \label{fig:IntersectionGamma}
\end{figure}

We consider the following flux function 
\begin{equation}
f(u) = \begin{cases}
d_1u + d_0, & x\leq u^* \text{ (Flux 1)}, \\
e_1u + e_0, & x > u^* \text{ (Flux 2)}
\end{cases}
\end{equation}
with the regularized flux function given by
\begin{equation}\label{eq:regflux}
f^\epsilon(u) = \begin{cases}
f_1(u) = d_1u + d_0, & 0 \leq u \leq u^*, \\
f^\epsilon_{\text{mid}}(u) = - \frac {1}{\epsilon}(f_1(u^*) - f_2^\epsilon(u^* + \epsilon))(u - u^*)+f_1(u^*), & u^*<u<u^*+\epsilon, \\
f_2^\epsilon(u) = e_1u+e_0, & u^*+\epsilon < u \leq u_{\text{max}}.
\end{cases}
\end{equation}
We define $u_\epsilon := u^* + \epsilon$. 
An example of the regularization is shown in figure \ref{fig:Linearization}.
Even though the regularization differs slightly from the one used in \cite{Wiensa}, the different cases to consider and the approach to solve the Riemann problems are completely analogous.
We refer to \cite[Section 3]{Wiensa} for further details.
We define $h: [0, u_\text{max}] \rightarrow \R$  to be a function connecting the points $(u^*, f(u^*+))$ and the initial value lying on Flux 2.
In addition, $\gamma$ denotes the intersection of $h$ and Flux 1,
see also figure \ref{fig:IntersectionGamma}.
The point $\gamma$ distinguishes the cases 3 and 4 in the following discussion. 
\subsubsection*{Case 1:  Either $u_L < u^*, u_R < u^*$ or $u_L > u^*, u_R > u^*$} This case corresponds to the classical case of solving Riemann problems, where the solution consists of a single rarefaction wave or shock, see e.g. \cite{Leveque2002}. 
\subsubsection*{Case 2: $ u_R < u^* < u_L$} By using the smallest convex hull approach the solution consists of a contact line following $f(u)$, connecting $u_R$ and $u^*$ and a shock, connecting $u^*$ and $u_L$.
The speed of the contact discontinuity is given by $d_1>0$ and the speed of the shock can be calculated with the Rankine-Hugoniot condition, i.e.,
\begin{align*}
s=\frac{f(u_L)-f(u^*)}{u_L-u^*}<0,
\end{align*}
where we recall that $f(u^*)=f(u^*-)$ holds.
The exact solution is then given by
\begin{equation}
u(x,t) = \begin{cases}
u_L, & \text{ if } x < s t, \\
u^*, & \text{ if }s t \leq x \leq d_1 t,\\
u_R, & \text{ if } x > d_1 t.
\end{cases}
\end{equation}
\subsubsection*{Case 3: $\gamma < u_L < u^* < u_R$}
Here, the solution is given by a shock connecting $u_L$ and $u^*$ and contact discontinuity reconnecting it to $u_R$.
The speed of the shock equals 
$$ s = \frac{f(u^*+)-f(u_L)}{u^*-u_L}<0.$$
Note that due to $u_L<u_R$ the convex hull construction chooses the flux $f(u^*+)$ at $u^*$.
The contact discontinuity connecting $u^*$ and $u_R$ moves at the speed $e_1<0$. 
Hence,
\begin{equation}
u(x,t) = \begin{cases}
u_L, & \text{ if } x < st, \\
u^*, & \text{ if }st \leq x \leq e_1 t, \\
u_R, & \text{ if } x >  e_1 t.
\end{cases}
\end{equation}
\subsubsection*{Case 4: $ u_L < u^* < u_R$ and $\gamma \geq u_L$}
In this case, we get only one shock connecting $u_L$ and $u_R$ as a solution due to the condition $\gamma \geq u_L$.
The speed of the shock can be calculated by the Rankin-Hugoniot condition.
The exact solution is given by,
\begin{equation}
u(x,t) = \begin{cases}
u_L, & \text{ if } x < st, \\
u_R, & \text{ if } x \geq st.
\end{cases}
\end{equation}

\begin{remark}
\rv{As aforementioned in \cite{Lu2009} Riemann solutions for piecewise quadratic discontinuous flux functions are derived. They also cover the case of a piecewise linear flux function if the quadratic terms are zero. For a general quadratic discontinuous flux, the solutions are more involved since no contact discontinuities occur.}
\end{remark}

Up to now, we have not addressed the case, where one of the boundary conditions equals the critical density $u^*$.
In this situation so called {\em zero waves}  can occur, see \cite{Wiensa}.
Those waves travel with speed $\mathcal{O}(\frac{1}{\epsilon})$ from right to left and have strength $\mathcal{O}(\epsilon)$, depending on the regularization parameter $\epsilon$.
Although they have no strength, information can be exchanged between Riemann problems lying next to each other influencing the solution.
For more details on the origin of zero waves in the context of traffic flow we refer to \cite{Wiensa}.
In particular, they arise when double Riemann problems are considered.
For a detailed discussion of these double Riemann problems we refer to \cite[section 4]{Wiensa}.

\section{LWR model with discontinuous flux on networks}
\label{sec:discnetworks}
Next, we focus on networks where we allow for  discontinuous flux functions.
The key idea is to consider the regularized flux function $f^\epsilon$ and take the limit $\epsilon \rightarrow 0$ to obtain results for the original problem at intersections.
In doing so, the standard approaches can be directly applied to the regularized flux function $f^\epsilon$. 

We start with a short introduction to the network setting.
For more details on traffic flow network models we refer the reader to \cite{garavellohanpiccoli2016book} and the references therein.

Let $G = (\mathcal{V}, \mathcal{E})$ be a directed finite graph consisting of a nonempty set $\mathcal{V}$ of vertices and a nonempty set of edges $\mathcal{E}$ representing the roads.
We denote by $\delta_v^- $ and $\delta_v^+$ the set of all incoming and outgoing edges for every vertex $v \in \mathcal{V}$.
Every edge $e \in \mathcal{E}$ is interpreted as an interval $I_e = (a_e,b_e) \subset \R$ representing the spatial extension of the road $e$.\\
For all edges $e$, which do not discharge into a vertex, i.e. $e \notin \cup_{v \in \mathcal{V}} \delta_v^-$, we set an infinite road length by $b_e = +\infty$.
Analogously, for all edges $e$, which do not originate from a vertex, i.e. $e \notin \cup_{v \in \mathcal{V}} \delta_v^+$, we set an infinite road length by $a_e = -\infty$. \\
We call the couple $(\mathcal{I}, \mathcal{V})$ with $\mathcal{I} = \{I_e: e\in \mathcal{E}\}$ a road network.

\rv{In order to derive the network solution, we restrict to the description of a single junction $v\in\mathcal{V}$. Hence, from now on we consider a fixed junction with $n$ incoming and $m$ outgoing roads.}
Then, we need additional information how vehicles are distributed from one road to another. 
Hence, we need to define a distribution matrix $ A \in [0,1]^{n \times m}$ as in \cite{garavellohanpiccoli2016book}: 
$$ A = \begin{pmatrix}
\beta_{1,1} & \cdots & \beta_{1,m}\\
\vdots & \vdots & \vdots \\
\beta_{n,1} & \cdots & \beta_{n,m}
\end{pmatrix}.$$
To conserve the mass we assume $\sum_{j = 1}^m \beta_{i,j} = 1$ for every $ i \in \{1, \dots, n\}$.
Furthermore, some technical assumptions on the distribution matrix are necessary which can be found in \cite[eq. (4.2.4)]{garavellohanpiccoli2016book}.
These assumptions are necessary to obtain unique solutions at the junction for the case $n\leq m$.
\rv{Following \cite[definition 4.2.4]{garavellohanpiccoli2016book} and definition \ref{WeakSolutionDisc},} weak solutions at a junction are intended in the following sense
\begin{definition}[Weak solution \rv{at a junction}]\label{WeakSolutionNetwork}
Let $v \in \mathcal{V}$ be a vertex with $n$ incoming and $m$ outgoing roads, represented \rv{by} the intervals $I_1^{\text{in}}, \dots I_n^{\text{in}}$ and $I_1^{\text{out}}, \dots I_n^{\text{out}}$. A weak solution at the junction $v$ is a collection $u$ \rv{and $w$} of functions $u_i^{\text{in}}: I_i^{\text{in}} \times \R_{\geq 0} \rightarrow \R$ \rv{, $w_i^{\text{in}}: I_i^{\text{in}} \times \R_{\geq 0} \rightarrow \R$ with $w_i^{\text{in}}(x,t)\in \tilde f(u_i^{\text{in}}(x,t))$ a.e.,} \mbox{where $i = 1, \dots,n$, $u_j^{\text{out}}$}: $I_j^{\text{out}} \times \R_{\geq 0} \rightarrow \R$ \rv{and $w_j^{\text{out}}: I_j^{\text{out}} \times \R_{\geq 0} \rightarrow \R$ with $w_j^{\text{out}}(x,t)\in \tilde f(u_j^{\text{out}}(x,t))$ a.e.,} where $ j = 1, \dots,m$, respectively, such that
\begin{align*} 
& \sum_{i = 1}^n \biggl( \int_0^T \int_{I_i^{\text{in}}} (u_i^{\text{in}}(x,t)\partial_t \phi_i^{\text{in}}(x,t) + \rv{w_i^{\text{in}}(x,t)}\partial_x\phi_i^{\text{in}}(x,t))\text{d}x\text{d}t \biggr)\\
+ &\sum_{j = 1}^m \biggl( \int_0^T \int_{I_j^{\text{out}}} (u_j^{\text{out}}(x,t)\partial_t \phi_j^{\text{out}}(x,t) + \rv{w_j^{\text{out}}(x,t)}\partial_x\phi_j^{\text{out}}(x,t))\text{d}x\text{d}t \biggr) = 0,
\end{align*}
for every collection of test functions $\phi_i^{\text{in}} \in C_0^1((a_i^{\text{in}},b_i^{\text{in}}] \times [0,T]), i = 1, \dots,n$ and $\phi_j^{\text{out}} \in C_0^1((a_j^{\text{out}},b_j^{\text{out}}] \times [0,T]), j = 1, \dots,m$, satisfying 
$$ \phi_i^{\text{in}}(b_i^{\text{in}}, \cdot) = \phi_j^{\text{out}}(a_j^{\text{out}}, \cdot), \quad \partial_x \phi_i^{\text{in}}(b_i^{\text{in}}, \cdot) = \partial_x \phi_j^{\text{out}}(a_j^{\text{out}}, \cdot), $$
for $ i \in {1,\dots n} \text{ and } j \in {1, \dots,m}.$
\end{definition}
Additionally, in order to get unique solutions, we will consider the following concept of admissible solutions, \rv{which adapts \cite[rule (A) and (B), \rv{p. 81}]{garavellohanpiccoli2016book} to the discontinuous setting:}
\begin{definition}[Admissbile Weak Solution]\label{AdmissibleWeakSolution}
We call $u$ \rv{and $w$ (as defined in definition \ref{WeakSolutionNetwork})} an admissible weak solution at a junction $v \in \mathcal{V}$ with the corresponding distribution matrix $A$ if 
\begin{enumerate}
\item $u_i^{\text{in}}(\cdot,t) \in \text{BV}(I_i^{\text{in}}), u_j^{\text{out}}(\cdot,t) \in \text{BV}(I_j^{\text{out}}) \text{ for all } i,j,$
\item $u$ \rv{and $w$ constitute} a weak solution at the junction $v$,
\item $\rv{w_j^{\text{out}}}(a_j^{\text{out}} +,t)) = \sum_{i = 1}^n \beta_{i,j} \rv{w_i^{\text{in}}}(b_i^{\text{in}} -,t)) \text{ for all } j = 1, \dots, m \text{ and } t \geq 0,$
\item $\sum_{i = 1}^n \rv{w_i^{\text{in}}}(b_i^{\text{in}} -,t)) $ is maximal subject to 2. and 3.
\end{enumerate}
\end{definition}
In particular, the maximization of the inflow with respect to the distributions parameters and the technical assumption \cite[eq. (4.2.4)]{garavellohanpiccoli2016book} guarantee the uniqueness of solutions for a continuous flux.

If $ n > m$, that means more incoming than outgoing roads, we need so-called \rv{right of way} parameters which prescribe in which order the vehicles are going to drive.
We consider the special case $n = 2$ and $m = 1$.
Let $C \in \rv{\R}$ denote the amount of vehicles that can move on the road. The \rv{right of way} parameter $q \in (0,1)$ describes the following situation:
If not all vehicles can move to the outgoing road, then $qC$ is the amount of vehicles that enter from the first and $(1-q)C$ the amount from the second road.
If there are more roads, we need to define several \rv{right of way} parameters.

For solving the maximization problems imposed by the definition \ref{AdmissibleWeakSolution} so-called supply and demand functions can be used, see \cite{lebacque1996godunov}.
The demand describes the maximal flux the incoming road wants to send.
In contrast, the supply describes the maximal flux the outgoing road is able to absorb.
The definition of the supply and demand functions of the regularized function is straightforward. 
As $f(u^*) \in \{ f(u^*+), f(u^*-)\}$ is not unique, we have to consider the limit for $\epsilon \rightarrow 0$ to determine the correct values of the supply and demand for the discontinuous flux function.
Nevertheless, it becomes apparent that the maximal flux in both cases is given by $f(u^*-)$.
This leads to the following definition:
\begin{definition}\label{def:supplydemand}
For a network with flux function $f$ having a discontinuity at $u^*$, the demand is given by
\begin{equation}
	D(u) = \begin{dcases*}
		f(u), &  $u \in [0, u^*),$\\
		f(u^*-), &  $u \in [u^*,u_{\text{max}}].$
		\end{dcases*}
\end{equation}

On the contrary, the supply reads as

\begin{equation}
	S(u) = \begin{dcases*}
f(u^*-), &  $u \in [0, u^*),$\\
f(u), &  $u \in (u^*,u_{\text{max}}]$\\
\end{dcases*}
\end{equation}
and
\begin{equation}
	S(u^*) = \begin{dcases*}
f(u^*-), &  \text{free flowing},\\
f(u^*+), &  \text{congested}.\\
\end{dcases*}
\end{equation}
\end{definition}
\begin{remark}
If we consider the regularized flux function $f^\epsilon$, the definition of supply and demand function is completely analogous by replacing $f$ with $f^\epsilon$.
\end{remark}

\begin{figure} 
  \centering
  \includegraphics[width=1\textwidth]{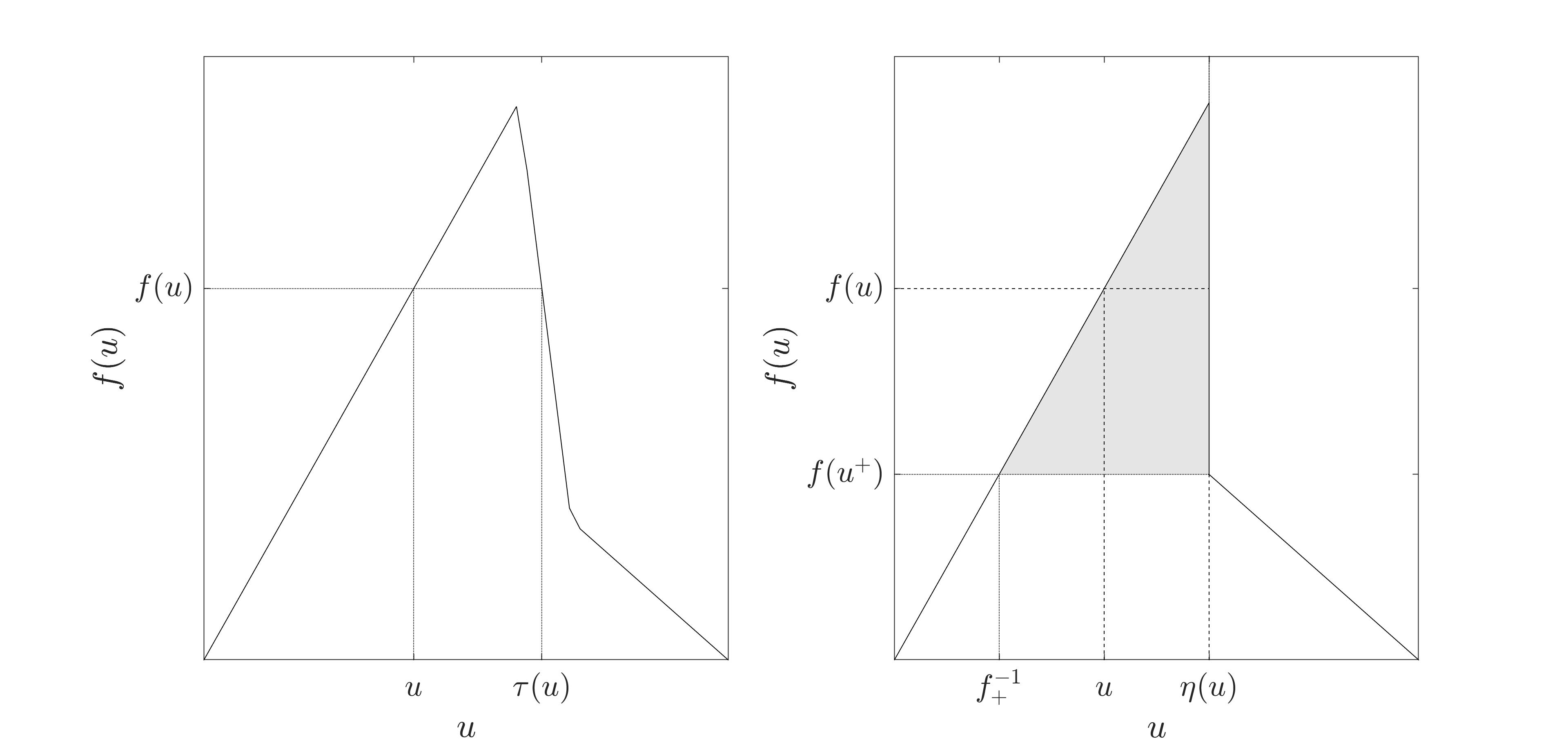}
   \caption[Function $\eta$]{The left panel shows the regularized version. The right panel shows the discontinuous flux $f$. The grey area displays the interval, where $f(\eta(u)) \leq f(u)$ and $\eta(u) = u^*$.}   
    \label{fig:eta}
\end{figure}

In order to show existence and uniqueness in the discontinuous case we need to define an additional function.
For a regularized flux function we notice that for every flux value, we get two different density values, see left picture in figure \ref{fig:eta}.
As different density values lead to different solutions, we need to be able to distinguish them and choose the correct solution.
\rv{In the continuous or regularized case a mapping usually called $\tau$ is used to determine the density values for which the flux values coincide, see also left picture in figure \ref{fig:eta} and \cite[definition 3.2.6]{garavellohanpiccoli2016book} for the precise definition of $\tau$.}
This situation becomes more complicated for the discontinuous flux functions.
For $\epsilon \rightarrow 0$ and $f(u) > f(u^*+)$, the corresponding density value \rv{given by $\tau$} can converge to the discontinuity area.
This situation is pictured in figure \ref{fig:eta} on the right.
\rv{Therefore, we need to adapt the definition of the mapping $\tau$ to also cover this case and determine the corresponding density value. This mapping is defined in the following and is denoted by $\eta$.}
\rv{For readability, we also define as an abbreviation $\fpustar:=f^{-1}(f(u^*+))$.}

\begin{definition}\label{eta_eps}
Let the function $\eta : [0, u_\text{max}] \rightarrow  [0, u_\text{max}]$ satisfy:
\begin{enumerate}
\item \rv{$f(\eta(u)) = f(u)$ for every $u \in \bigl([0, \rv{\fpustar}] \cup [u^*,u_\text{max}]\bigr) $ and $ \eta(u) \neq u $ for every $u \in \bigl([0, \rv{\fpustar}] \cup (u^*,u_\text{max}]\bigr) $.}  
\item For $u \in \bigl(\rv{\fpustar}, u^*\bigr)$ it holds, $\eta(u) = u^*$ and $f(\eta(u))=f(u)$.
\end{enumerate}
\end{definition}

\noindent Note that for the case $u \in \bigl(\rv{\fpustar}, u^*\bigr)$ the flux value is given by $f(u)\in(f(u^*+),f(u^*-))$. 

\begin{remark}
\rv{
We note that the mapping $\tau$ from  \cite[definition 3.2.6]{garavellohanpiccoli2016book} is in principle given by the first statement in definition \ref{eta_eps} with the separated interval $[0, \rv{\fpustar}] \cup [u^*,u_\text{max}]$ being replaced by $[0,u_\text{max}]$.
Hence, the only difference between $\tau$ and $\eta$ lies in the treatment of $u \in \bigl(\rv{\fpustar}, u^*\bigr)$.
}
\end{remark}

\subsection{Riemann solutions at a junction}

Now, we present a Riemann solver for two types of junctions.
First, we consider a junction with $n$ incoming and $m$ outgoing roads with $n\leq m$, second, a junction with $n = 2$ incoming and $m = 1$ outgoing roads.
For the continuous case the existence of solutions satisfying in particular the definition \ref{AdmissibleWeakSolution} has been shown in \cite{garavellohanpiccoli2016book}.
\rv{Since the continuous flux is not multivalued, the definition \ref{AdmissibleWeakSolution} coincides with \cite[Definition 4.2.4]{garavellohanpiccoli2016book}.}
As we have seen in section \ref{sec_ZeroWaveRP}, the solution to a Riemann problem may produce more than one wave.
In order to have a well-defined network problem it is important that the waves produced at the junction have negative speed on the incoming and positive speed on the outgoing roads.
This is ensured by the following theorem:

\begin{theorem}\label{thm:ntom}
Let $w \in \mathcal{V}$ be a junction with $n$ incoming and $m$ outgoing roads \rv{with $n\leq m$}. For every constant initial datum $u_{i,0}^{\text{in}} , u_{j,0}^{\text{out}}  \in [0, u_{\text{max}}],\ i=1,\ldots,n,\ j=1,\ldots,m$, there exists a unique admissible weak solution $u$ \rv{in the sense of definition \ref{AdmissibleWeakSolution}.}
Moreover, for the incoming road the solution is given by the wave $(u_{i,0}^{\text{in}} , u_{i}^{\text{in}}),\ i=1,\ldots,n$ and for every outgoing road the solution is the wave $(u_{j,0}^{\text{out}},u_{j}^{\text{out}}), j=1,\ldots,m$, where

\begin{equation}\label{eq:ntomuin}
u_i^{\text{in}} \in \begin{dcases*}
\{ u_{i,0}^{\text{in}} \} \cup (\eta(u_{i,0}^{\text{in}}),u_{\text{max}}], &  if $u_{i,0}^\text{in} \in [0, \rv{\fpustar}]$, \\
\{u_{i,0}^{\text{in}}\} \cup [u^*,u_{\text{max}}],  &if $ u_{i,0}^\text{in} \in (\rv{\fpustar}, u^*)$, \\
[u^*,u_{\text{max}}], & if $u_{i,0}^\text{in} \in [u^*, u_\text{max}],$
\end{dcases*}
\end{equation}
and 
\begin{equation}\label{eq:ntomuout}
u_{j}^{\text{out}} \in \begin{dcases*}
[0, u^*], & if $u_{j,0}^\text{out} \in [0,u^*)$,\\
[0,u^*],& if  $u_{j,0}^{\text{out}}=u^*$ \text{and free flowing},\\
\{u^*\}\cup [0, \rv{\fpustar}),& if $u_{j,0}^{\text{out}}=u^*$ \text{and congested},\\
\{ u_{j,0}^{\text{out}} \} \cup [0, \eta(u_{j,0}^{\text{out}})), &  if $u_{j,0}^\text{out} \in (u^*,u_{\max}].$
\end{dcases*}
\end{equation}
\end{theorem}
\begin{proof}
Using the definition of the supply and demand functions in definition \ref{def:supplydemand} and \rv{the results from \cite[section 5.2.3]{garavello2006traffic} we can follow the proof of \cite[theorem 5.1.2]{garavello2006traffic}} and uniquely determine the inflows which maximize the flux through junctions subject to the distribution parameters. It remains to show that by the choice of the density values the correct waves are induced. 

We start with considering the outgoing roads.
If $u_{j,0}^{\text{out}}\in [0,u^*)$, the choice $u_j^{\text{out}}\in[0,u^*]$ produces the classical case of a shock or contact discontinuity, both with the positive speed $d_1$, since the flux is a linear increasing function. 
If $u_{j,0}^{\text{out}}=u^*$, we need to know if the traffic ahead is congested or not in order to determine the correct flux value of $u^*$ and the corresponding density.
In particular, we are now concerned with a double Riemann problem by $u_j^{\text{out}}$, $u_{j,0}^{\text{out}}=u^*$ and the information about the traffic ahead. By the analysis done in \cite[section 4]{Wiensa}
the choices of $u_j^{\text{out}}$ give rise to waves with positive velocities, even though the velocities differ depending on the traffic situation. Note that the free flowing case corresponds to case 1 in \cite[section 4]{Wiensa} and the congested case to case 2. In particular, in the latter case $f_j^{\text{out}}\leq f(u^*+)=f(u^*)$ holds such that the solution is a shock with velocity always greater or equal to zero, compare also case 4 of section \ref{sec_ZeroWaveRP}.

For $u_{j,0}^{\text{out}}\in (u^*,u_{\text{max}}]$ we need to choose the density $u_{j}^{\text{out}} \in \{ u_{j,0}^{\text{out}} \} \cup [0, \eta(u_{j,0}^{\text{out}}))$ which corresponds to the incoming flux.
Here, we are in case 4 of section \ref{sec_ZeroWaveRP}, as $u_{j,0}^{\text{out}}<\gamma$, such that the solution is a shock with positive velocity.

On the contrary, considering the incoming roads and $u_{i,0}^\text{in} \in [0, \rv{\fpustar}]$, the choice $u_i^{\text{in}}\in \{ u_{i,0}^{\text{in}} \} \cup (\eta(u_{i,0}^{\text{in}}),u_{\text{max}}]$ also corresponds to case 4, but this time the shock has a negative velocity.

Now, let $u_{i,0}^\text{in} \in (u^*, u_\text{max}]$ and $f_i^{\text{in}}$ denote the uniquely determined inflow.
If $f_i^{\text{in}}<f(u^*+)$, the choice $u_i^{\text{in}}\in (u^*, u_\text{max}]$ corresponds to the solution of a Riemann problem with a linear decreasing function and hence either a shock or contact discontinuity with negative velocity.
On the contrary if $f_i^{\text{in}}\geq f(u^*+)$, we choose $u_i^{\text{in}}=u^*$.
Here, the convex hull reconstruction induces a shock as a solution with velocity
\begin{align*}
s=\frac{f_i^{\text{in}}-f(u_{i,0}^{\text{in}})}{u^*-u_{i,0}^{\text{in}}}<0.
\end{align*}
Further, if $u_{i,0}^\text{in}=u^*$, we would normally also need to consider a double Riemann problem for which we would need the information of boundary data of the incoming roads. Nevertheless, the solutions at the junction are the same. As can be seen in \cite[section 4]{Wiensa} cases 2 and 4, if $f_i^{\text{in}}<f(u^*+)$ the choice $u_i^{\text{in}}\in (u^*, u_\text{max}]$ induces a contact discontinuity with negative velocity as the solution at the junction.
If $f_i^{\text{in}}\geq f(u^*+)$ we choose $u_i^{\text{in}}=u^*$ which gives a constant solution.

Now, let us turn to the remaining case of $u_{i,0}^\text{in} \in (\rv{\fpustar}, u^*)$.
Again the solution depends on the inflow.
If $f_i^{\text{in}}<f(u^*+)$, the choice of $u_{i}^{\text{in}}\in\{u_{i,0}^{\text{in}}\} \cup (u^*,u_{\text{max}}]$ induces negative waves.
Nevertheless, depending on $u_{i,0}^{\text{in}}> \gamma$ or $u_{i,0}^{\text{in}}\leq \gamma$ we are either in case 3 or case 4 of section \ref{sec_ZeroWaveRP}.
Even though case 3 induces two waves, all have a negative velocity.
If $f_i^{\text{in}}\geq f(u^*+)$, choosing $u_i^{\text{in}}=u^*$ gives a shock with velocity
\begin{align*}
s=\frac{f_i^{\text{in}}-f(u_{i,0}^{\text{in}})}{u^*-u_{i,0}^{\text{in}}}\leq 0,
\end{align*}
as $f_i^{\text{in}}$ is bounded from above by $f(u_{i,0}^{\text{in}})$ due to the definition of the demand.

Hence, the choices of the densities induce the correct waves.
\end{proof}

Now, we consider the case of more incoming than outgoing roads.
Exemplary, we study the 2-to-1 situation, even though the results can be easily extended to the $n$-to-1 case.
We assume that a \rv{right of way} parameter $q$ is given and the unique density of the discontinuous flux can be determined according to the following theorem: 

\begin{theorem}\label{thm:2to1}
Let $w \in \mathcal{V}$ be a junction with $n = 2$ incoming and $m = 1$ outgoing roads with a right way parameter $q \in (0, 1)$. For every constant initial datum $u_{i,0}^{\text{in}} , u_{1,0}^{\text{out}} \in [0, u_{\text{max}}]$ with $i \in \{1,2\}$, there exists a unique admissible weak solution $u$ \rv{in the sense of definition \ref{AdmissibleWeakSolution}.}
Moreover, for every incoming road $i$ the solution is given by the wave $(u_{i,0}^{\text{in}} , u_{i}^{\text{in}})$ and for the outgoing edge the solution is the wave $(u_{1,0}^{\text{out}},u_{1}^{\text{out}})$, where

\begin{equation}\label{eq:2to1uin}
u_i^{\text{in}} \in \begin{dcases*}
\{ u_{i,0}^{\text{in}} \} \cup (\eta(u_{i,0}^{\text{in}}),u_{\text{max}}], &  if $u_{i,0}^\text{in} \in [0, \rv{\fpustar}]$, \\
\{u_{i,0}^{\text{in}}\} \cup [u^*,u_{\text{max}}],  &if $ u_{i,0}^\text{in} \in (\rv{\fpustar}, u^*)$, \\
[u^*,u_{\text{max}}], & if $u_{i,0}^\text{in} \in [u^*, u_\text{max}],$
\end{dcases*}
\end{equation}
and

\begin{equation}\label{eq:2to1uout}
u_{1}^{\text{out}} \in \begin{dcases*}
[0, u^*], & if $u_{1,0}^\text{out} \in [0,u^*)$,\\
[0,u^*],& if  $u_{1,0}^{\text{out}}=u^*$ and free flowing,\\
\{u^*\}\cup [0, \rv{\fpustar}),& if $u_{1,0}^{\text{out}}=u^*$ and congested,\\
\{ u_{1,0}^{\text{out}} \} \cup [0, \eta(u_{1,0}^{\text{out}})), &  if $u_{1,0}^\text{out} \in (u^*,u_{\max}].$
\end{dcases*}
\end{equation}
\end{theorem}

\begin{proof}
Following \cite[Section 3.2.2]{garavellohanpiccoli2016book}, the flux values at the junction can be calculated with the following steps:
\begin{enumerate}
\item Calculate the maximal possible flux $f_\text{max} = \min\{ D(u_{1,0}^\text{in}) + D(u_{2,0}^\text{in}), S(u_{1,0}^\text{out})\}$.
\item Consider the \rv{right of way} parameter and the flux maximization and calculate the intersection $P= (qf_\text{max}, (1-q)f_\text{max})$. 
\item If $P$ is in the feasible area $\Omega =  \{ z\in [0, D(u_{1,0}^{in})] \times [0, D(u_{2,0}^\text{in})]: z_1 + z_2 \leq f_\text{max}\}$, set $(z_1,z_2)=P$.
If not, determine the point in $\Omega \cap \{(z_1,z_2): z_1+z_2=f_{\text{max}}\}$ which is the closest to the line $z_2=(1-q)/q z_1$.
\item The flux values are given by $f^\text{in}_1 = z_1,\ f^\text{in}_2 = z_2,\ f^\text{out}_1 = z_1+z_2$.
\end{enumerate}

Completely analogous to theorem \ref{thm:ntom} we can show that the choice of the densities admits the correct wave speeds.
\end{proof}
The Riemann solutions proposed in theorem \ref{thm:ntom} and theorem \ref{thm:2to1} are the key ingredients for the splitting algorithm on networks in the next section. 

\section{A splitting algorithm on networks}
\label{sec:splitnetworks}
Different problems might occur when designing a numerical scheme for a conservation law with  discontinuous flux.
However, the main difficulties are induced by the zero waves.
Since these waves have infinite speed, the regular CFL condition is scaled by the regularization parameter $\epsilon$.
This leads to very inefficient step sizes and high computational effort.
Possible ways to avoid using the regular CFL condition are implicit methods \cite{ImplicitMethod}
or the splitting algorithm introduced in \cite{Towers}, see also \cite{BurgerChalonsVilladaMulticlass}.
In this section we recall the main ideas of the splitting algorithm and also present its extension to networks.

\begin{figure}[htb]
  \centering
  \includegraphics[width=1\textwidth]{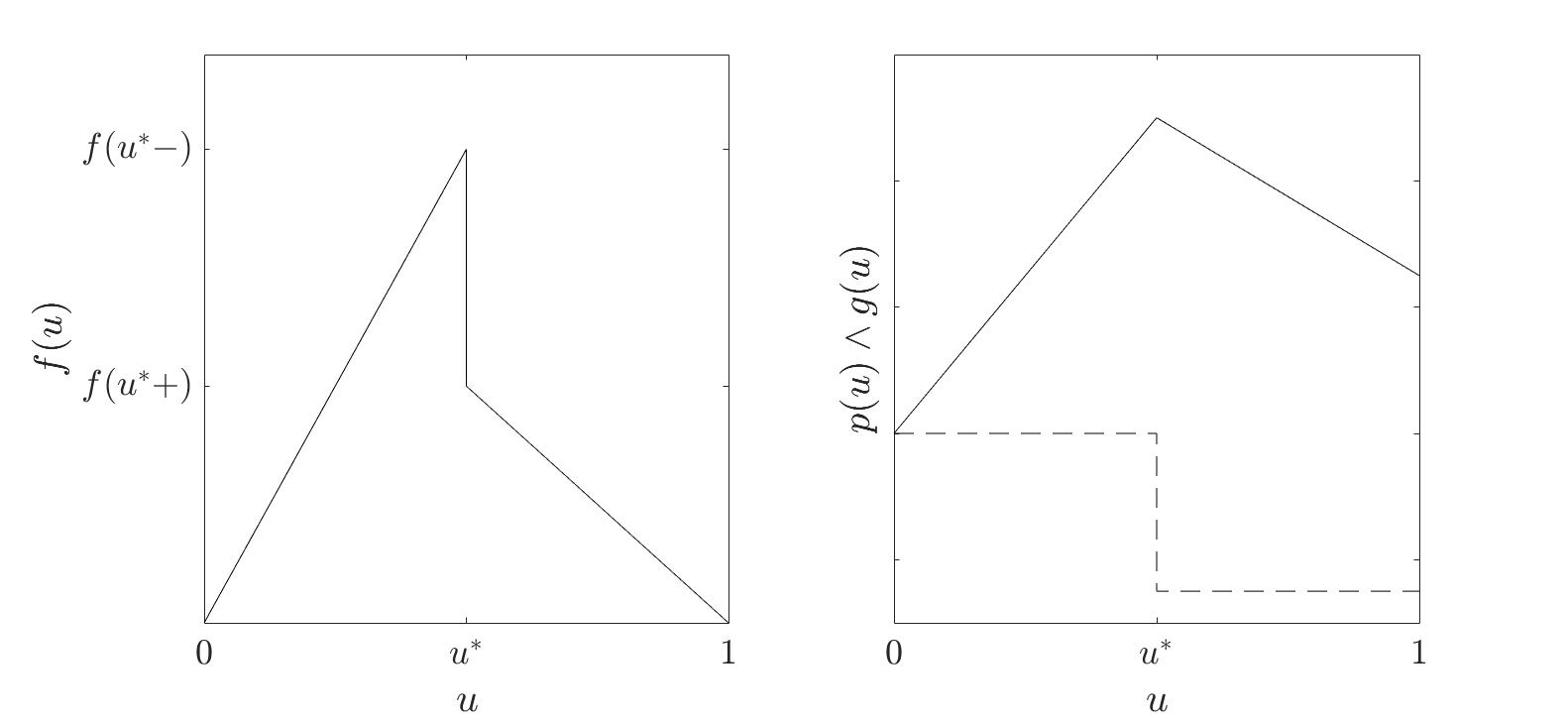}
   \caption[Splitting Algorithm]{Based on \cite[figure 1]{Towers}. The left panel shows the discontinuous flux function, the right panel shows the splitted version. The dotted line represents the piecewise constant function $g$, the solid line shows the continuous function $p$ .}   
    \label{fig:DiscFluxSplitting}
\end{figure}

We consider a flux function $f$ with discontinuity at $u^*$, which satisfies assumption \ref{assumptionSplit}.
The jump has magnitude $\alpha = f(u^*-)-f(u^*+)$.
The main idea of the splitting algorithm is to split the discontinuous function $f$ into two parts, namely $p$ and $g$, such that $f(u) = p(u) + g(u)$.
Here, $g$ is a piecewise constant function defined by
\begin{equation*}
g(u) = -\alpha H(u-u^*),
\end{equation*}
where $H$ denotes the Heaviside function.
This function extracts the jump from the function $f$. 
In addition, we define $p(u) = f(u) - g(u).$
Hence, $p$ is continuous and can be used for numerical algorithms.
An example can be seen in figure \ref{fig:DiscFluxSplitting}.

\subsection{A single road}

This case has been already treated in \cite{Towers}
and will be the basis for the splitting algorithm on networks.
When solving the scalar conservation law \eqref{ScalarConservationLaw} on a single road, the boundary value in the case $s(t) = u^*$ is not unique.
In the continuous case, the flux value is determined by $\mathcal{F}(t)$ in \eqref{eq:fluxboundary}.
When splitting the function $f$, the discontinuity is shifted to the function $g$.
So analogously, we define the multivalued version of $g$ in order to be able to assign a flux value in the critical case, i.e.,
\begin{equation}
\tilde g(u) = \begin{cases}
0, & u \in [0, u^*),\\
[-\alpha, 0], & u = u^*, \\
-\alpha, & u \in (u^*, u_{\text{max}}].
\end{cases}
\end{equation}
Furthermore, we define $\mathcal{G}(t) \in \tilde g(s(t)).$
It holds, 
\begin{equation} \mathcal{F} (t) = f(u^*-) \Leftrightarrow \mathcal{G}(t) = 0, \qquad \mathcal{F}(t) = f(u^*+) \Leftrightarrow \mathcal{G}(t) = -\alpha.
\end{equation}
Additionally to the assumptions \ref{assumptionSplit}, we assume:
\begin{assumption}\cite[Assumption 1.1]{Towers}\label{assumptionSplit2}
The initial data satisfies $u_0(x) \in [0, u_{\text{max}}]$ for $x \in (a,b),\ u_0 \in \text{BV}([a,b])$ and $g(u_0) \in \text{BV}([a,b]).$ 
We also assume that $\mathcal{G}(t) \in [-\alpha,0]$ for $t \in [0,T]$, and $\mathcal{G} \in \text{BV}([0,T])$.
\end{assumption}
\begin{remark} \rv{We emphasize that the original splitting algorithm for a single road \cite{Towers} is not limited to piecewise linear discontinuous flux functions. Another prominent example might be concave piecewise quadratic flux functions with discontinuity again at $u^*$, cf. \cite{Lu2009,Towers,Wiensa}}.

\end{remark}
Then, we are able to handle the flux function $p$ and the jump $g$ separately. 
As usual, we discretize the domain $\Pi_T$ by a spatial mesh and a temporal mesh. 
Let $\Delta x$ denote the spatial and $\Delta t$ the temporal step size.
We define the grid constant via $$\lambda = \frac{\Delta t}{\Delta x}.$$ 
For an integer $\rv{K} \in \N$ the grid size is given by $\Delta x = \frac{b-a}{\rv{K}+1}$ and the grid points by $x_{\rv{k}} = a+\rv{k} \Delta x$.
Further, we define the $\rv{K}$ disjunct intervals $I_{\rv{k}} = [x_{\rv{k}} - \frac{\Delta x}{2}, x_{\rv{k}} + \frac{\Delta x}{2})$ and $\mathcal{\rv{K}} = \{ 1, \cdots \rv{K}\}$ and $\mathcal{\rv{K}}^+ = \{ 0, \cdots \rv{K}+1\}$.\\
The time steps are defined by $t^n = n \Delta t$, for $n = 0,1, \dots,N$. The value $N \in \N$ is chosen such that the time horizon fulfills $T \in [t^N, t^N+\Delta t).$ 

As the algorithm splits the function $f$ into two parts, $f(u) = p(u) +  g(u)$, every step of the algorithm consists of two half steps.
Here, the first half step corresponds to approximately solving $ u_t + g(u)_x = 0$ denoted by $ U_{\rv{k}}^{n+\frac{1}{2}} \approx u(x_{\rv{k}},t^n).$
Consequently, the second step approximates $u_t + p(u)_x = 0$ and is denoted by $ U_{\rv{k}}^n  \approx u(x_{\rv{k}},t^n).$

We denote the backward spatial difference by $\Delta_- U_{\rv{k}}^n = U_{\rv{k}}^n - U_{\rv{k}-1}^n$ and the forward spatial difference by $\Delta_+ U_{\rv{k}}^n = U_{\rv{k}+1}^n - U_{\rv{k}}^n$.
The initial values are denoted by $U_{\rv{k}}^0 = u_0(x_{\rv{k}})$, the boundary values by
\begin{align*}
r^n &= r^{n+\frac{1}{2}} = r(t^n),  &U_0^{n+\frac{1}{2}} &= U_0^n = r^n \\
s^n &= s^{n+ \frac{1}{2}} = s(t^n) & U_{\rv{K}+1}^{n+\frac{1}{2}} &= U_{\rv{K}+1}^n = s^n.
\end{align*}

The function $g$ covers the issue of $s(t^n) = u^*$.
We define
\begin{align}\label{boundgonstreet}
g_{\rv{K}+1}^{n+\frac{1}{2}} = g_{\rv{K}+1}^n = \mathcal{G}(t^n) = \begin{cases}
0, & \text{ if } s(t^n) < u^*, \\
0, & \text{ if } s(t^n) = u^*, \text{ traffic ahead of } x = b \text{ is free-flowing}, \\
-\alpha, & \text{ if }  s(t^n) = u^*, \text{ traffic ahead of } x = b \text{ is congested}, \\
-\alpha, & \text{ if } s(t^n) > u^*.
\end{cases}
\end{align}
That means, we can describe the boundary value $\mathcal{F}(t)$ via $\mathcal{G}(t)$.
For the splitting algorithm in its original fashion \cite{Towers} we additionally need the following function.

\begin{definition}\cite[Eq. (3.7)]{Towers}
Let $G(z) := z - \lambda g(z)$ with grid constant $\lambda$ and $\tilde G$ its multivalued version. The function $\tilde G$ is strictly increasing and has a unique inverse, $z \rightarrow \tilde G^{-1}(z)$, which is a single valued function. It is Lipschitz continuous and nondecreasing. It holds, $0 \leq \tilde G^{-1}(z) \leq z$. The functions are given by

\begin{align*}
\tilde G(u) = 
\begin{cases}
u, & u \in [0,u^*), \\
[u^*, u^*+\lambda \alpha], & u = u^*, \\
u + \lambda \alpha, & u \in (u^*,u_\text{max}],
\end{cases}
\qquad
\tilde G^{-1} (u) =
\begin{cases}
u, & u \in [0, u^*),\\
u^*, &u \in [u^*,u^*+\lambda \alpha), \\
u-\lambda \alpha, & u \in [u^*+\lambda \alpha, u_\text{max} + \lambda \alpha].
\end{cases}
\end{align*}
\end{definition}

The splitting algorithm \cite{Towers} can be then expressed as
\begin{equation}\label{GleichungSplittAlgScheme}
\begin{cases}

\begin{cases}
U_{\rv{k}}^{n+1/2} = \tilde G^{-1}\bigl(U_{\rv{k}}^n - \lambda g_{\rv{k}+1}^{n+1/2}\bigr), & \qquad \qquad  \rv{k} = \rv{K},\rv{K}-1, \ldots, 1, \\
g_{\rv{k}}^{n+1/2} = \bigl(U_{\rv{k}}^{n+1/2} - U_{\rv{k}}^n + \lambda g_{\rv{k}+1}^{n+1/2}\bigl)/\lambda, & \qquad \qquad  \rv{k = \rv{K},\rv{K}-1, \ldots, 1}, \\
\end{cases} \\[3ex]
U_{\rv{k}}^{n+1} = U_{\rv{k}}^{n+1/2} - \lambda \Delta_{-}p^g\bigl(U_{\rv{k}+1}^{n+1/2},U_{\rv{k}}^{n+1/2}\bigr), \qquad \quad \ \rv{k} \in \mathcal{\rv{K}}.
\end{cases}
\end{equation}
Note that the first half step, which includes the first two equations, is implicit.
Nevertheless, instead of solving a nonlinear system of equations, the equation can be solved backwards in space starting with $\rv{k} = \rv{K}$.
The last equation in (\ref{GleichungSplittAlgScheme}) uses the Godunov scheme for the second half step which flux function is denoted by $p^g$.

We note that for the implicit equation a CFL condition is not needed, but it is required for the third step.
As $p$ is continuous, we can use the standard CFL condition of the Godunov scheme. Note that similar to \cite{Towers} also other two-point monotone schemes, e.g. the Lax-Friedrichs method, can be used.

As shown in \cite[Theorem 5.1]{Towers} the splitting algorithm \eqref{GleichungSplittAlgScheme} converges to a weak solution of \eqref{ScalarConservationLaw}.
However, obtaining a similar statement about weak entropy solutions is still an open problem.

\subsection{Networks}
The key idea to numerically solve such discontinuous conservation laws on networks is to use the splitting algorithm only on the roads and determine the correct in- and outflows at the boundaries by the help of the Riemann solver established in, e.g. theorem \ref{thm:ntom}.
As the splitting algorithm works with flux values, there is no need to compute the exact densities at the junction.
Instead we need to know how the solution at the junction influences the flux values.
\rv{The algorithm for a single junction is depicted in algorithm \ref{algo:complete}.
The general description of the algorithm allows for either junctions with given distribution or right of way parameters.
For simplicity, we assume that each road is represented by the same interval $(a,b)$, such that we have $I_1^{\text{in}}=\dots=I_n^{\text{in}}=I_1^\text{out}=\dots=I_m^\text{out}=(a,b)$.
\begin{remark}
Note that this simplification enables the use of the same grid points on each road which spares further sub- or superindices. However, the algorithm can be easily adapted to different road lengths.
\end{remark}
\noindent We assume in the following that the space and time grid is the same as in the previous subsection.
\begin{algorithm}[htb!]
	\caption{Splitting algorithm for a fixed junction of dispersing or merging type}\label{algo:complete}
	\begin{algorithmic}[1]
		\REQUIRE number of incoming roads $n$, number of outgoing roads $m$, distribution matrix $A$, (right of way parameters), flux $f$, discontinuity $u^*$, jump magnitude $\alpha$, interval boundaries $a,\ b$, number of grid points $K$, time step size $\Delta t$, initial values $u^\text{in}_{i,0}(x)$ and $u^\text{out}_{j,0}(x)$, boundary values $u^\text{in}_{i,0}(a,t)$ and $u^\text{out}_{j,0}(b,t)$
		\ENSURE approximate solutions $U_{i,k}^{n,\text{in}}$ and $U_{j,k}^{n,\text{out}}$
		\STATE Initilization:
		\STATE $\Delta x=(b-a)/(K+1)$ and $\lambda=\Delta t/\Delta x$,
		\STATE $U_{i,k}^{0,\text{in}}=u^\text{in}_{i,0}(x_k)$ for $i=1,\dots,m$, $U_{j,k}^{0,\text{out}}=u^\text{out}_{j,0}(x_k)$ for $j=1,\dots,n$
		\STATE $U_{i,0}^{n,\text{in}}=u^\text{in}_{i,0}(a,t^n)$, $U_{j,0}^{n,\text{out}}=u^\text{out}_{j,0}(b,t^n)$
		\FOR{$n=0,\dots,N-1$}
		\STATE Solve the by definition \ref{AdmissibleWeakSolution} induced optimization problem at the junction based on the flux $f$, which gives the fluxes $f_i^\text{in},\ f_j^{\text{out}}$
		\STATE Compute the densities at the junction with an appropriate Riemann solver
		\STATE Compute the adjusted flux values for the incoming roads $f_{i,\text{adj.}}^\text{in}$
		\FOR{$i=1,\dots,n$}
		\STATE $g_{i,K+1}^{n+1/2,\text{in}} = g_{i,K+1}^{n,\text{in}} =f_i^{\text{in}}-f_{i,\text{adj.}}^\text{in}$
			\FOR{$k=K,K-1,\dots,1$}
			\STATE $U_{i,k}^{n+1/2,\text{in}} = \tilde G^{-1}\bigl(U_{i,k}^{n,\text{in}} - \lambda g_{i,k+1}^{n+1/2,\text{in}}\bigr)$
			\STATE $g_{i,k}^{n+1/2,\text{in}} = \bigl(U_{i,k}^{n+1/2,\text{in}} - U_{i,k}^{n,\text{in}} + \lambda g_{i,k+1}^{n+1/2,\text{in}}\bigl)/\lambda$
			\ENDFOR
					\STATE $p^g\bigl(U_{i,K+1}^{n+1/2,\text{in}},U_{i,K}^{n+1/2,\text{in}}\bigr)=f_{i,\text{adj.}}^\text{in}$
			\FOR{$k\in \mathcal{K}$}
			\STATE $U_{i,k}^{n+1,\text{in}} = U_{i,k}^{n+1/2,\text{in}} - \lambda \Delta_{-}p^g\bigl(U_{i,k+1}^{n+1/2,\text{in}},U_{i,k}^{n+1/2,\text{in}}\bigr)$
			\ENDFOR
			\ENDFOR
		\FOR{$j=1,\dots,m$}
		\STATE Compute $g_{j,K+1}^{n+1/2,\text{out}} = g_{j,K+1}^{n,\text{out}}$ as in \eqref{boundgonstreet}
			\FOR{$k=K,K-1,\dots,1$}
			\STATE $U_{j,k}^{n+1/2,\text{out}} = \tilde G^{-1}\bigl(U_{j,k}^{n,\text{in}} - \lambda g_{j,k+1}^{n+1/2,\text{out}}\bigr)$
			\STATE $g_{j,k}^{n+1/2,\text{out}} = \bigl(U_{j,k}^{n+1/2,\text{out}} - U_{j,k}^{n,\text{out}} + \lambda g_{j,k+1}^{n+1/2,\text{out}}\bigl)/\lambda$
			\ENDFOR
			\STATE $p^g\bigl(U_{j,1}^{n+1/2,\text{in}},U_{j,0}^{n+1/2,\text{out}}\bigr)=f_j^\text{out}$
			\FOR{$k\in \mathcal{K}$}
			\STATE $U_{j,k}^{n+1,\text{out}} = U_{j,k}^{n+1/2,\text{out}} - \lambda \Delta_{-}p^g\bigl(U_{j,k+1}^{n+1/2,\text{out}},U_{j,k}^{n+1/2,\text{out}}\bigr)$
			\ENDFOR
			\ENDFOR
		\ENDFOR
	\end{algorithmic}
\end{algorithm}
The approximate solutions are  denoted by $U_{i,k}^{n,\text{in}}\approx u_{i,k}^{n,\text{in}}(x_k,t^n)$ for $i=1,\dots,n$ and $U_{j,k}^{n,\text{out}}\approx u_{j,k}^{n,\text{out}}(x_{\rv{k}},t^n)$ for $i=j,\dots,m$.
The half steps are denoted accordingly.
We assume that the initial values $u^\text{in}_{i,0}(x)$ and $u^\text{out}_{j,0}(x)$ and the boundary values $u^\text{in}_{i,0}(a,t)$ and $u^\text{out}_{j,0}(b,t)$ are given.
Then, we can directly set the starting values as $U_{i,k}^{0,\text{in}}=u^\text{in}_{i,0}(x_k)$ for $i=1,\dots,n$ and $U_{j,k}^{0,\text{out}}=u^\text{out}_{j,0}(x_k)$ for $j=1,\dots,m$.
The left boundary values for the incoming roads are given by $U_{i,0}^{n,\text{in}}=u^\text{in}_{i,0}(a,t^n)$ and the right boundary values for the outgoing roads by $U_{j,0}^{n,\text{out}}=u^\text{out}_{j,0}(b,t^n)$.
The missing boundary data, i.e. the right one for the incoming road and the left one for the outgoing road, are determined by the flux values at the junction.
\newline The overall strategy of the splitting algorithm on a network consists of three important steps:
\begin{enumerate}
\item Solve the optimization problem induced by definition \ref{AdmissibleWeakSolution} (in particular item 4) at the junction to calculate the flux values $f_i^\text{in}, f_j^\text{out}$ with the discontinuous flux $f$.
\end{enumerate}
Here, it is crucial to use the discontinuous flux function $f$ when calculating demand and supply at the junction, which are necessary for the optimization problem, since otherwise not only the flux values can be different but also the decision whether supply or demand is active. Therefore, in line 6, we need to compute first the fluxes at the junction with the original discontinuous flux function $f$.
\newline These flux values bring us now to:
\begin{itemize}
\item[2.] Consider the corresponding density values at the junction applying \rv{either theorem \ref{thm:ntom} or theorem \ref{thm:2to1} depending on the type of junction}. If the density value is greater or equal than $u^*$, i.e. the road is congested, the flux value needs to be adjusted. Hence, we determine adjusted flux values $f_\text{adj.}^\text{in/out}$.
\end{itemize}
Using the calculated (unadjusted) flux values from step one, we can determine the densities at the junction with the help of the appropriate Riemann solver (theorem \ref{thm:ntom} and \eqref{eq:ntomuin}-\eqref{eq:ntomuout} or theorem \ref{thm:2to1} and \eqref{eq:2to1uin}-\eqref{eq:2to1uout}) at the junction (line 7).
Then, these density values can be used to calculate the corresponding flux value of $p$ to get the boundary values and decide about an adjustment.
This means that if the density given by the Riemann solver suggests that the traffic ahead is free flowing, no adjustment needs to be made - as the flux value already lies on $p$.
If the traffic ahead is congested we take the flux value on the curve $p$ for the corresponding density value and not its original flux value on the curve $f$.
\newline Both the first step, i.e solving the optimization problem at the junction, and the second step, i.e. how to adjust the flux, strongly depend on the junction type.
Hence, we postpone a detailed discussion to the following subsections, in which we will consider the 1-to-1, 1-to-2 and 2-to-1 case in more detail.
\newline In addition, the first steps give us all the ingredients for the final step:
\begin{itemize}
\item[3.] Away from the junction, the splitting algorithm can be used. At the junction itself, the adjusted fluxes are used to determine the missing boundary data. 
\end{itemize}
The adjusted flux values from the previous step are important for the second half step (line 17 or 28) of the splitting algorithm which uses a Godunov type scheme based on $p$.
Here, the missing fluxes in and out of the junction the Godunov scheme needs are simply given by these adjusted flux values.
The flux values can be used directly as boundary conditions, see lines 15 and 26, and the exact density values are not needed for the scheme.
\newline Further boundary data is needed in the first half step of the algorithm, lines 12-13 and 23-24.
Here, we start with $k = K$ to avoid solving a nonlinear equation system.
For example, the value $g^{n,\text{out}}_{j,K+1}$ is determined by the density value at the boundary.
Nevertheless, for every outgoing road the boundary values are known and the  value $g^{n,\text{out}}_{j,K+1}$ can be determined as before in \eqref{boundgonstreet}.
However, for the incoming road, i.e. for $g^{n,\text{in}}_{i,K+1}$, the density value at the junction is necessary.
From the adjustment made in step two we know whether the traffic ahead is congested or not.
In addition, instead of working with the exact density values at the junction directly, we can also use the adjusted flows computed in the second step (or line 8).
Here, we can simply compute the boundary value $g^{n,\text{in}}_{i,K+1}$ by the difference between the fluxes calculated with $f$ (before adjustment) and the ones with $p$ (after adjustment),  i.e.
\begin{align}\label{eq:gnetwork}
g^{n+1/2,\text{in}}_{i,K+1} = g^{n,\text{in}}_{i,K+1} =f_i^{\text{in}}-f_{i,\text{adj.}}^\text{in}
\end{align}
Note that the definition of $g^{n,\text{in}}_{i,K+1}$ in \eqref{eq:gnetwork} is slightly different in comparison to \eqref{boundgonstreet} or for the outgoing roads.
Nevertheless, if the density value on the incoming road is smaller than $u^*$, $g^{n,\text{in}}_{i,K+1}$ equals zero (no adjustment needs to be made) and if it is greater $u^*$, it equals $-\alpha$.
Only in the case the density value equals $u^*$ the situation is more involved, as the traffic on the incoming road can be congested and $f(u^*-)>f^{\text{in}}>f(u^*+)$ can occur such that we need $g^{n,\text{in}}_{i,K+1}\in (-\alpha,0)$.
\newline Furthermore, we can decrease the computational costs of the algorithm: In the second step (or in line 7 of algorithm \ref{algo:complete}) the density at the junction is computed.
This can be very expensive and hence we aim to avoid this.
In the third step we have seen that for the missing boundary data only the adjusted flux values are necessary and not the densities at the junction themselves.
Therefore, studying first each junction type in detail allows to determine the corresponding flux values based on the density values and the supply and demand functions and the intermediate expensive step for the computation of the exact densities can be skipped.
Hence, we combine the first and second step of the strategy in one single step.
In the following, we will study a 1-to-1 junction in detail and present the tailored algorithm.
As the strategy is completely analogous for the 1-to-2 junction and 2-to-1 junction, we will only present the algorithms and discuss important properties.
The algorithms can then be used to replace the lines 6 to 8 in algorithm \ref{algo:complete}.
\begin{remark}
The extension to $n$-to-1 or 1-to-$m$ junctions follows immediately from the upcoming discussion.
For more complex junctions, it is possible to adapt the proposed strategy by studying the different solutions of the optimization problem at the junction.
\end{remark}
\begin{remark}
Further note that the presented strategy and also algorithm \ref{algo:complete} can be used for arbitrary junctions and nonlinear discontinuous flux functions once an appropriate Riemann solver similar to theorem \ref{thm:ntom} and \ref{thm:2to1} is established.
\end{remark}
}
\subsubsection{1-to-1 junction}
First, we consider a junction with one incoming and one outgoing road in detail.
Let $u_{1,0}^\text{in}, u_{1,0}^\text{out}$ be constant initial values.
The flux values are given by $f^{\text{in}}_1=f^{\text{out}}_1=\min\{D(u_{1,0}^\text{in}),S(u_{1,0}^\text{out})\}$.
We distinguish the following cases:
\subsubsection*{Case A: demand and supply are equal }
There are two different situations depending on the density value of the incoming road where demand and supply can be equal.
\begin{enumerate}
\item \textbf{ $\mathbf{u_{1,0}^\text{in}} \mathbf{\geq u^*}$: } Here, the traffic on the outgoing road needs to be free-flowing or equal to $u^*$ and the traffic ahead of the outgoing road is free-flowing.
Supply and demand equal the maximal possible flux $f(u^*-)$.
For these flux values the solution of the Riemann solver on the incoming and outgoing road is given by $u^*$.
As corresponding flux values already lie on $p$, we do not need to adjust them.
\item \textbf{ $\mathbf{u_{1,0}^\text{in}} \mathbf{< u^*}$: }  Demand and supply can also be equal, if the number of vehicles the demand wants to send is equal to the number of vehicles the supply can take.
Therefore, $u_{1,0}^\text{in} < \gamma$ and $u_{1,0}^{\text{out}}>u^*$ need to hold.
The solution proposed by the Riemann solver at the junction is the initial condition for each road.
Hence, the traffic on the outgoing road stays congested and we need to adjust the flux value on the outgoing road, i.e., $f_{1,\text{adj}}^\text{out} = f_1^\text{out} + \alpha$.
The flux on the incoming road stays free-flowing.
\end{enumerate}

\subsubsection*{Case B: supply is restrictive}
If the supply is restrictive, i.e. $f^{\text{in}}_1=f^{\text{out}}_1=S(u_{1,0}^\text{out})<D(u_{1,0}^\text{in})$, the outgoing road cannot take the whole amount of cars, the demand wants to send.
From this follows that the traffic on the outgoing road is congested, i.e. $u_{1,0}^\text{out} > u^*$ or  $u_{1,0}^\text{out} = u^*$ and the traffic ahead \rv{being congested}.
This has an effect on the incoming road.
Again, we face two situations:
\begin{enumerate}
\item \textbf{ $\mathbf{u_{1,0}}^\text{in} \mathbf{< u^*}$: }   Initially, the traffic on the incoming road is free-flowing and the traffic on the outgoing road is congested.
As the supply is restrictive, the traffic on the incoming road congests as well, the solution of the Riemann solver is on both roads given by $u_{1,0}^\text{out}$.
Therefore, we have to adjust the flux value, $f_{1,\text{adj}}^\text{in} = f_1^\text{in}+\alpha.$
In addition, we need to set $\rv{g_{1,K+1}^{n,\text{in}}} = -\alpha$.
The outgoing flux needs to be also adjusted, i.e., $f_{1,\text{adj}}^\text{out} = f_1^\text{out}+\alpha.$

\item \textbf{ $\mathbf{u_{1,0}}^\text{in} \mathbf{> u^*}$: } If the traffic on the incoming road is congested as well, the solution of the Riemann solver on both roads is again given by $u_{1,0}^\text{out}$.
This leads to the same adjustments as in the previous case.
\end{enumerate}
\subsubsection*{Case C: demand is restrictive}
If the demand is restrictive, the outgoing road is able to take the whole amount of vehicles the demand wants to send. 
Here, we only have one possible situation for the initial condition on the incoming road:
\begin{enumerate}
\item \textbf{ $\mathbf{u_{1,0}}^\text{in} \mathbf{\leq u^*}$: }  
The solution of the Riemann problem on each road is now given by $u_{1,0}^\text{in}$, such that no adjustment is necessary.
\end{enumerate}
Note that we have the same flux function on each road.
Hence, in the case ${u_{1,0}}^\text{in} {> u^*}$ either the supply is restrictive or demand and supply are equal.

The whole procedure is summarized in algorithm \ref{algo:1to1}.

\begin{algorithm}[htb]
	\caption{1-to-1 junction}\label{algo:1to1}
	\label{Flux1_1_in_D_S}
	\begin{algorithmic}[1]
		\REQUIRE  Demand $D_1$, Supply $S_1$, flux $f$, discontinuity $u^*$, jump magnitude $\alpha$
		\ENSURE Flux values $f_1^\text{in}, f_1^\text{out}$
		\STATE $f_1^\text{in} = \min \bigl\{D_1, S_1\bigr \}$
		\STATE  $f_1^\text{out} = f_1^\text{in} $
		\IF{$S_1 = D_1$ \textbf{and} $D_1 \neq f(u^*-)$}						
					\STATE $f_1^\text{out} = f_1^\text{out} + \alpha$
					\ELSIF {$S_1<D_1$}
					\STATE {$f_1^\text{out} = f_1^\text{out} + \alpha$}
					\STATE {$f_1^\text{in} = f_1^\text{in} + \alpha$}
		\ENDIF
	\end{algorithmic}
\end{algorithm}

\begin{remark}
\rv{Theoretically, the Riemann solver in theorem \ref{thm:ntom} coincides for a 1-to-1 junction with the Riemann solver on a single road. Hence, the procedure described in algorithm \ref{algo:1to1} leads to the same solution.
In contrast to that on the numerical level, the splitting algorithm for a 1-to-1 junction does not exactly coincide with the splitting algorithm used for a single road \cite{Towers}. 
The reason for the computational difference is that the splitting algorithm for the 1-to-1 junction considers the exact solution of the Riemann problem at the junction point and hence uses exact values for $g_{1,K+1}^{n,\text{in}}$ while for a single road this value is only approximated using $g_{1,K+1}^{n,\text{out}}$.}
\end{remark}
\subsubsection{1-to-2 junction}
The difference to the 1-to-1 junction is now that we have to consider two supplies. So the case distinctions to determine the flux values are a bit more complex. However, the procedure itself does not change.
The details can be seen in algorithm \ref{algo:1to2}.

\begin{algorithm}[h!]
	\caption{1-to-2 junction}
	\label{algo:1to2}
	\begin{algorithmic}[1]
		\REQUIRE  Demand $D_1$, Supply $S_1, S_2$, distribution parameter $\beta_1, \beta_2$, flux $f$, discontinuity $u^*$, jump magnitude $\alpha$
		\ENSURE Flux values $f_1^\text{in}, f_1^\text{out}, f_2^\text{out}$
		\STATE $f_1^\text{in} = \min \bigl\{D_1, S_1/\beta_1, S_2/\beta_2 \bigr \}$
		\STATE  $f_1^\text{out} = \beta_1 f_1^\text{in}, f_2^\text{out} = \beta_2 f_1^\text{in} $
		\IF {$f_1^\text{in} = D_1$}
		\IF{$S_1 = D_1$ \textbf{and} $S_1 \neq f(u^*-)$}						
					\STATE $f_1^\text{out} = f_1^\text{out} + \alpha$
		\ELSIF{$S_2 = D_1$ \textbf{and} $S_2 \neq f(u^*-)$}						
					\STATE $f_2^\text{out} = f_2^\text{out} + \alpha$
		\ENDIF
		\ELSIF {$f_1^\text{in}=S_1/\beta_1$ \textbf{or} $f_1^\text{in}=S_2/\beta_2$}
		\IF {$f_1^\text{in}>f(u^*+)$}
		\STATE {$f_1^\text{in}=f(u^*-)$}
		\ELSE
		\STATE {$f_1^\text{in}=f_1^\text{in}+\alpha$}
		\ENDIF
		\IF {$S_1/\beta_1=S_2\beta_2$ \textbf{and} $S_1\beta_1<D_1$}
		\STATE $f_1^\text{out} = f_1^\text{out} + \alpha$
		\STATE $f_2^\text{out} = f_2^\text{out} + \alpha$
		\ELSIF {$S_1/\beta_1<S_2\beta_2$ \textbf{and} $S_1\beta_1<D_1$}
		\STATE $f_1^\text{out} = f_1^\text{out} + \alpha$
		\ELSIF {$S_2/\beta_2<S_1\beta_1$ \textbf{and} $S_2\beta_2<D_1$}
		\STATE $f_2^\text{out} = f_2^\text{out} + \alpha$
		\ENDIF
		\ENDIF
	\end{algorithmic}
\end{algorithm}

We remark that if the inflow on the first road equals the demand and the restriction given by at least one of the supplies, the latter road needs to be congested.
The Riemann solver states congestion such that the flux needs to be adjusted. Then, the incoming road either stays free flowing or the solution is given by $u^*$.
Here, no adjustment is necessary.
Note that the conditions $S_i\neq f(u^*-)$ are necessary for the specific case of $\beta_i=1$.

If at least one of the supply restrictions is active, we need to adjust the corresponding flux values as in the 1-to-1 case and the incoming road.
Nevertheless, here an interesting case (which is not possible in the 1-to-1 situation) can occur.
The solution of the Riemann problem at the incoming road can be given by $u^*$ and the inflow is in $f^{\text{in}}_1\in(f(u^*+),f(u^*-))$, see also the proof of theorem \ref{thm:ntom}.
We need to adjust the flux to $p(u^*)=f(u^*-)$.
Hence, in this case if we use \eqref{eq:gnetwork} the adjustment is between $0$ and $-\alpha$.
Finally, as in the 1-to-1 case there is no need for an adjustment if the demand is restrictive.
 
%
\begin{algorithm}[h!]
	\caption{2-to-1 junction}
	\label{algo:2to1}
	\begin{algorithmic}[1]
		\REQUIRE  Demand $D_1$, $D_2$, Supply $S_1$, right of way parameter $q_1, q_2$, flux $f$, discontinuity $u^*$, jump magnitude $\alpha$
		\ENSURE Flux values $f_1^\text{in}, f_2^\text{in}, f_1^\text{out}$
		\STATE $f_\text{max} = \min \bigl\{D_1+D_2, S_1 \bigr \}$
		\STATE  $f_1^\text{out}= f_\text{max}, z_1= q_1 f_\text{max}, z_2=q_2 f_\text{max} $
		\IF {$z_1>D_1$}
			\STATE {$z_1=D_1, z_2=f_\text{max}-z_2$}
		\ELSIF {$z_2>D_2$}
			\STATE {$z_2=D_2, z_1=f_\text{max}-z_1$}
		\ENDIF
		\STATE $f_1^\text{in}=z_1, f_2^\text{in}=z_2$		
		\IF {$D_1+D_2 = S_1$ \textbf{and} $S_1 \neq f(u^*-)$}
			\STATE $f_1^\text{out} = f_1^\text{out} + \alpha$
		\ELSIF {$D_1+D_2<S_1$}
			\IF{$S_1 \neq f(u^*-)$}						
					\STATE $f_1^\text{out} = f_1^\text{out} + \alpha$
			\ENDIF
			\IF {$f_1^\text{in}<D_1$}
				\IF {$f_1^\text{in}>f(u^*+)$}
				\STATE {$f_1^\text{in}=f(u^*-)$}
				\ELSE
				\STATE {$f_1^\text{in}=f_1^\text{in}+\alpha$}
			\ENDIF
\ENDIF
		\IF {$f_2^\text{in}<D_2$}
			\IF {$f_2^\text{in}>f(u^*+)$}
			\STATE {$f_2^\text{in}=f(u^*-)$}
			\ELSE
			\STATE {$f_2^\text{in}=f_2^\text{in}+\alpha$}
			\ENDIF
		\ENDIF			
		\ENDIF
	\end{algorithmic}
\end{algorithm}

\subsubsection{2-to-1 junction}
Recall that for two incoming and one outgoing roads, the maximal possible flux on the outgoing road is given by $f_\text{max} = \min\{ D(u_{1,0}^\text{in}) + D(u_{2,0}^\text{in}), S(u_{1,0}^\text{out})\}$.
The corresponding flux values on the incoming roads are always smaller or equal to the actual demand.
The algorithm for the corresponding adjusted flux values is listed in algorithm~\ref{algo:2to1}.

As before, no adjustment is needed when the demand on both roads is smaller than the supply.
If the supply is active we might need to adjust the outgoing road and in most cases at least one incoming road. 
As in the 1-to-2 case, $f_i^{\text{in}}$ can be greater than $f(u^*+)$ such that the Riemann solver of theorem \ref{thm:2to1} gives $u^*$ as a solution with the corresponding flux $p(u^*)=f(u^*-)$ leading to an adjustment smaller than $\alpha$.

\section{Numerical Simulation}
\label{sec:num}
In this section, we present some numerical examples to compare the splitting algorithm on networks with the Riemann solver.
Further, we compute the solution using a regularized flux and the Godunov scheme.
We consider the following flux function
\begin{align}
f(u)=\begin{cases}
u,& u\in [0, 0.5),\\
0.5(1-u),& u\in[0.5,1].
\end{cases}
\end{align}
The corresponding regularization is given by \eqref{eq:regflux}.

We consider in particular the 1-to-2 and 2-to-1 situations.
For our comparison, we choose constant initial data on each road.
The junction is always located at $x=0$, such that the incoming roads have negative $x$-coordinates and the outgoing ones positive.
All roads are assumed to have the same length of 2.
Furthermore, we compute the sum of the $L^1$ errors on each road at $t=1$ to compare the numerical approaches.
As the splitting algorithm requires $\lambda\leq 1$ as CFL condition, but the regularization approach $\lambda\leq \epsilon$, we compare different CFL constants and regularization parameters with each other.

\subsection{1-to-2 junction}
In both scenarios the supply of the first outgoing road is restrictive.
The parameter settings are as follows:
\begin{enumerate}
\item[1)] $u_{1,0}^\text{in} = 0.4$, $u_{1,0}^\text{out} = 0.9$, $u_{2,0}^\text{out} = 0.7$,  $\beta = (0.75, 0.25)$. \\
The exact solution is given by
\begin{align*}
u_1^\text{in}(x,t) = \begin{cases}
0.4, & x\geq -\frac{3}{2}t, \\
0.5, & -\frac{3}{2}t<x<-\frac{1}{2}t,\\
\frac{13}{15}, & -\frac{1}{2}t < x <0,
\end{cases}
\qquad 
&u_1^\text{out}(x,t) = 
0.9, \qquad  \\
&u_2^\text{out}(x,t) = \begin{cases}
\frac{1}{15}\beta_2  & 0 \leq x\leq \frac{8}{41}t, \\
0.7, & \frac{8}{41}t<x. 
\end{cases}
\end{align*}
Apparently, the solution induces two waves on the incoming road.
\item[2)] $u_{1,0}^\text{in} = 0.4$, $u_{1,0}^\text{out} = 0.7$, $u_{2,0}^\text{out} = 0.2$,  $\beta = (0.5, 0.5)$. \\
The exact solution is given by
\begin{align*}
u_1^\text{in}(x,t) = \begin{cases}
0.4, & x\geq -t, \\
0.5, & -t<x<0,
\end{cases}
\qquad 
&u_1^\text{out}(x,t) = 
0.7, \qquad 1 \leq x \leq 3, \\
&u_2^\text{out}(x,t) = \begin{cases}
0.3\beta_2  & 1 \leq x\leq t, \\
0.2, & t<x. 
\end{cases}
\end{align*}
This example generates $f_1^{\text{in}}=0.3$ such that we are in the specific case that the flux has to be adjusted, but only to $0.5$.
\end{enumerate}

In table \ref{Convergence_Net_Ex1}, the $L^1$-errors are compared for different $\Delta x$ and the numerical convergence rate (CR) is determined by a least square fitting method.
For the splitting algorithm we are allowed to choose a rather large CFL constant of 0.75 while for a direct comparison to the regularized approach $\lambda=0.1$. is chosen.
Furthermore, the regularization approach is also computed with a smaller regularization parameter.
For the regularization approach we choose $\lambda=\epsilon$.

\begin{table}[htb]
\centering
\begin{tabular}{ | cx{1.5pt} c | c | c | c | }
\hline
 &\multicolumn{4}{c|}{Example 1 } \\
 \hline
 & \multicolumn{2}{c|}{Splitt.} &\multicolumn{2}{c|}{Reg.}\\
 \hline
	$\Delta x$	& $\lambda=0.75$&$\lambda=0.1$ & $\epsilon = 0.1$ 			& $\epsilon = 0.01$		 	\\
 \hline
$0.04$& 33.44e-03 & 46.77e-03 & 82.26e-03& 51.82e-03\\
$0.02$& 24.17e-03 &29.05e-03 &65.59e-03 & 30.69e-03\\
$0.01$& 14.16e-03 &20.12e-03 & 60.86e-03& 24.45e-03\\
$0.005$& 8.97e-03 & 12.49e-03& 58.37e-03& 20.44e-03\\
\hline
CR & 0.64695 & 0.62453 &0.1593 &0.4353 \\
\hline
\hline
& \multicolumn{4}{c|}{Example 2 }\\
\hline
 & \multicolumn{2}{c|}{Splitt.} &\multicolumn{2}{c|}{Reg.}\\
 \hline
	$\Delta x$ & $\lambda=0.75$&$\lambda=0.1$ & $\epsilon = 0.1$ 			& $\epsilon = 0.01$	\\
\hline
$0.04$ & 4.58e-03& 7.41e-03 &44.70e-03 &14.57e-03\\
$0.02$ & 2.97e-03& 4.24e-03&43.47e-03&11.28e-03\\
$0.01$ & 2.03e-03&2.89e-03 &42.50e-03&10.06e-03\\
$0.005$ &  1.24e-03&1.99e-03 &41.80e-03&9.21e-03\\
\hline
CR &0.61911& 0.62327 &0.0322&0.2150\\
\hline
\end{tabular}
\begin{tikzpicture}[overlay]
\draw[ultra thick] (-6.1,0.15) rectangle (-2.1,3);
\draw[ultra thick] (-6.1,-0.95) rectangle (-2.1,-3.8);
\end{tikzpicture}
\caption{The table shows the $L^1$ error for the splitting algorithm (Splitt.) and the regularized problem (Reg.) for two different $\epsilon$ values for the 1-to-2 examples.}
\label{Convergence_Net_Ex1}
\end{table}

We can see that the error terms obtained by
the splitting algorithm are the lowest and so the computational costs.
Obviously, the error terms increase with a lower CFL due to numerical diffusion.
For a direct comparison with the regularized approach, the CFL condition should be the same.
\rv{Meaning that the second column in table \ref{Convergence_Net_Ex1} for the splitting algorithm should be compared with the first one of the regularized approach.}
In this case, we can see that the splitting algorithm also performs better in both examples.
By choosing a smaller regularization parameter the performance of the regularized approach increases, but also the computational costs.
To obtain similar error terms as for the splitting algorithm the regularization parameter needs to be further reduced at very high computational costs.

\begin{figure}[htb]
  \centering
  \includegraphics[width=1\textwidth]{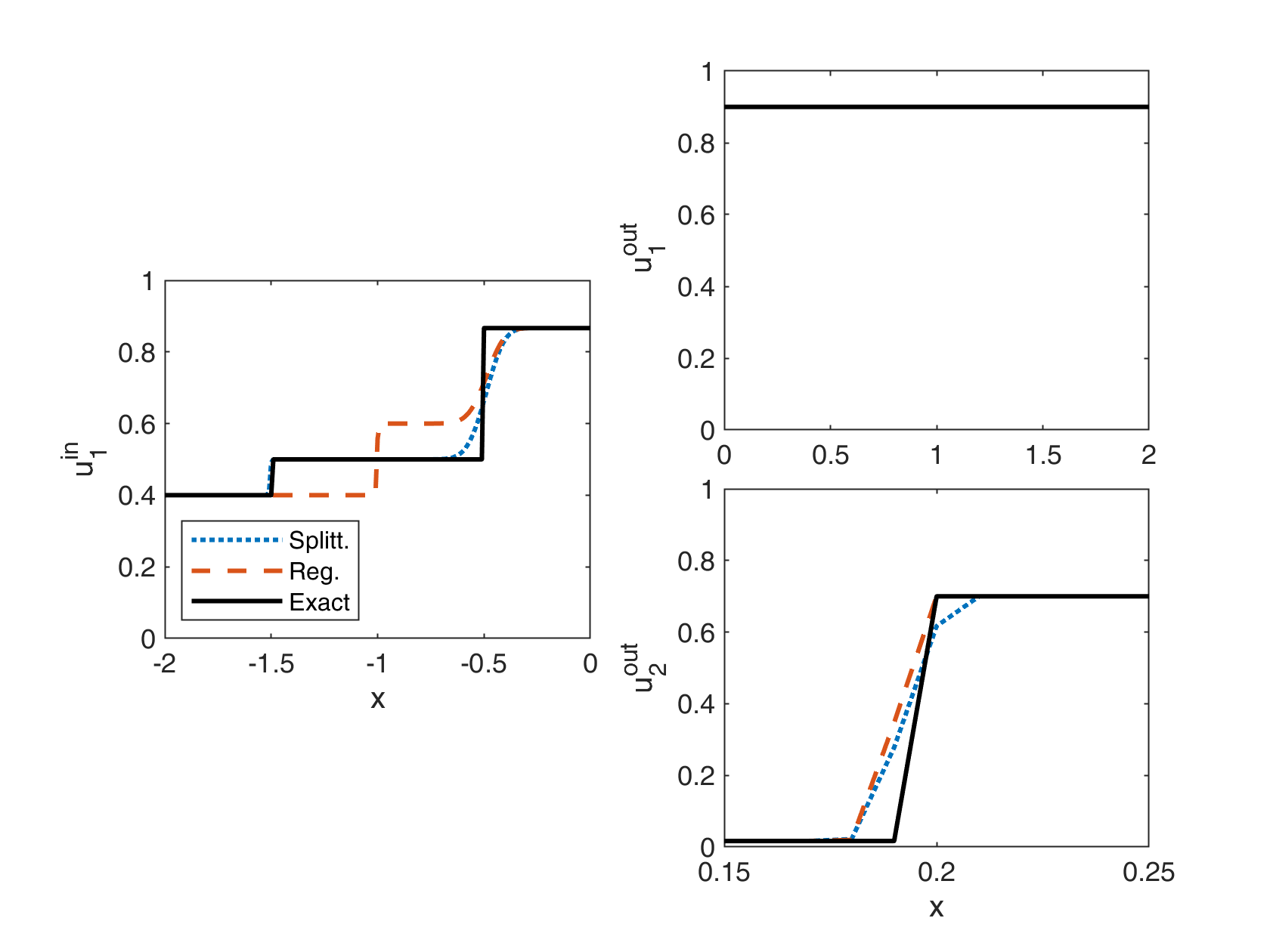}
   \caption{First example with $\lambda=0.1$ and $\Delta x=0.01$.}   
    \label{fig:1to2}
\end{figure}

For  $\lambda=0.1$ and $\Delta x=0.01$ the approximate solutions and the exact solution are displayed in figure \ref{fig:1to2} for the first example.
On the incoming road the splitting algorithm approximates the solution much better, while for the regularized approach an even smaller regularization parameter would be needed for a correct approximation.
On the outgoing roads both schemes capture the right waves.
Note that we zoomed in on the second outgoing road to make the small differences visible.

\subsection{2-to-1 junction}
Here, we consider two scenarios, where in the first scenario the demand is restrictive while in the second one the supply. 
The parameter settings are as follows:
\begin{enumerate}
\item[1)] $u_{1,0}^\text{in} = 0.2$, $u_{2,0}^\text{in} = 0.25$, $u_{1,0}^\text{out} = 0.3$,  $q = 0.75$. \\
The exact solution is given by
\begin{align*}
&u_1^\text{in}(x,t) =0.2,
\qquad 
&u_1^\text{out}(x,t) = \begin{cases}
0.45  & 0 \leq x\leq t, \\
0.3, & t<x
\end{cases}\\
&u_2^\text{in}(x,t)=0.25. 
& 
\end{align*}
\item[2)] $u_{1,0}^\text{in} = 0.6$, $u_{2,0}^\text{in} = 0.7$, $u_{1,0}^\text{out} = 0.4$,  $q = 0.8$. \\
The exact solution is given by
\begin{align*}
&u_1^\text{in}(x,t) = \begin{cases}
0.5  & x\leq -2t, \\
0.6, & -2t<x<0
\end{cases},
\qquad 
&u_1^\text{out}(x,t) = \begin{cases}
0.5  & 0 \leq x\leq t, \\
0.4, & t<x
\end{cases},\\
&u_2^\text{in}(x,t)=\begin{cases}
0.8  & x\leq -0.5 t, \\
0.7, & -0.5 t<x<0
\end{cases}.
& 
\end{align*}
In particular, the flux value for the first incoming road needs to be adjusted from $0.4$ to $0.5$ in the splitting algorithm, while on the second road it is adjusted as usual from $0.1$ to $0.35=0.1+\alpha$.
Note that due to the high velocity  on the first incoming road we evaluate the error at $t=0.5$. 
\end{enumerate}

\begin{table}[htb]
\centering
\begin{tabular}{| cx{1.5pt}c| c | c| c |}
\hline
 &\multicolumn{4}{c|}{Example 1 } \\
 \hline
 & \multicolumn{2}{c|}{Splitt.} &\multicolumn{2}{c|}{Reg.} \\
 \hline
	$\Delta x$	& $\lambda=0.75$&$\lambda=0.1$ & $\epsilon = 0.1$ 			& $\epsilon = 0.01$		 	\\
 \hline
$0.04$& 9.25e-03 & 16.22e-03 & 16.22e-03&17.01e-03 \\
$0.02$& 5.90e-03 &11.63e-03 &11.63e-03 & 12.19e-03\\
$0.01$& 2.98e-03 &8.13e-03 & 8.13e-03& 8.52e-03\\
$0.005$& 8.97e-03 &  5.71e-03& 5.71e-03& 5.99e-03\\

\hline

CR & 0.53838 & 0.50353 &0.50353 &0.50347 \\
\hline
\hline
 & \multicolumn{4}{c|}{Example 2 }\\
 \hline
 & \multicolumn{2}{c|}{Splitt.} &\multicolumn{2}{c|}{Reg.}\\
 \hline
$\Delta x$ & $\lambda=0.75$&$\lambda=0.1$ & $\epsilon = 0.1$ 			& $\epsilon = 0.01$	\\
\hline
$0.04$& 14.12e-03& 20.10e-03 &85.69e-03 &27.55e-03\\
$0.02$& 9.65e-03& 13.86e-03&79.98e-03&21.65e-03\\
$0.01$& 6.41e-03&9.57e-03 &75.96e-03&17.49e-03\\
$0.005$& 4.51e-03&6.69e-03 &73.22e-03&14.65e-03\\
\hline
CR &0.55295& 0.52959 &0.07551&0.30432\\
\hline
\end{tabular}
\begin{tikzpicture}[overlay]
\draw[ultra thick] (-6.1,0.15) rectangle (-2.1,3);
\draw[ultra thick] (-6.1,-0.95) rectangle (-2.1,-3.8);
\end{tikzpicture}
\caption{The table shows the $L^1$ error for the splitting algorithm (Splitt.) and the regularized problem (Reg.) for two different $\epsilon$ values for the 2-to-1 examples.} 
\label{Convergence_Net_Ex2}
\end{table}

Considering the $L^1$ errors in table \ref{Convergence_Net_Ex2} similar observations as in the 1-to-2 case can be made.
The splitting algorithm performs best and has the lowest computational costs.
In the first example the solution of the splitting algorithm coincides with the regularized approach as the flows are equal at the junction due to the active demand, \rv{see second column of the splitting algorithm and first one of the regularized approach in table \ref{Convergence_Net_Ex2}}.
Nevertheless, in the second example the splitting algorithm performs much better for the same CFL constant.
In addition, the regularization parameter needs to be even smaller than 0.01 to reach the precision of the faster splitting algorithm. 

\begin{figure}[htb]
  \centering
  \includegraphics[width=1\textwidth]{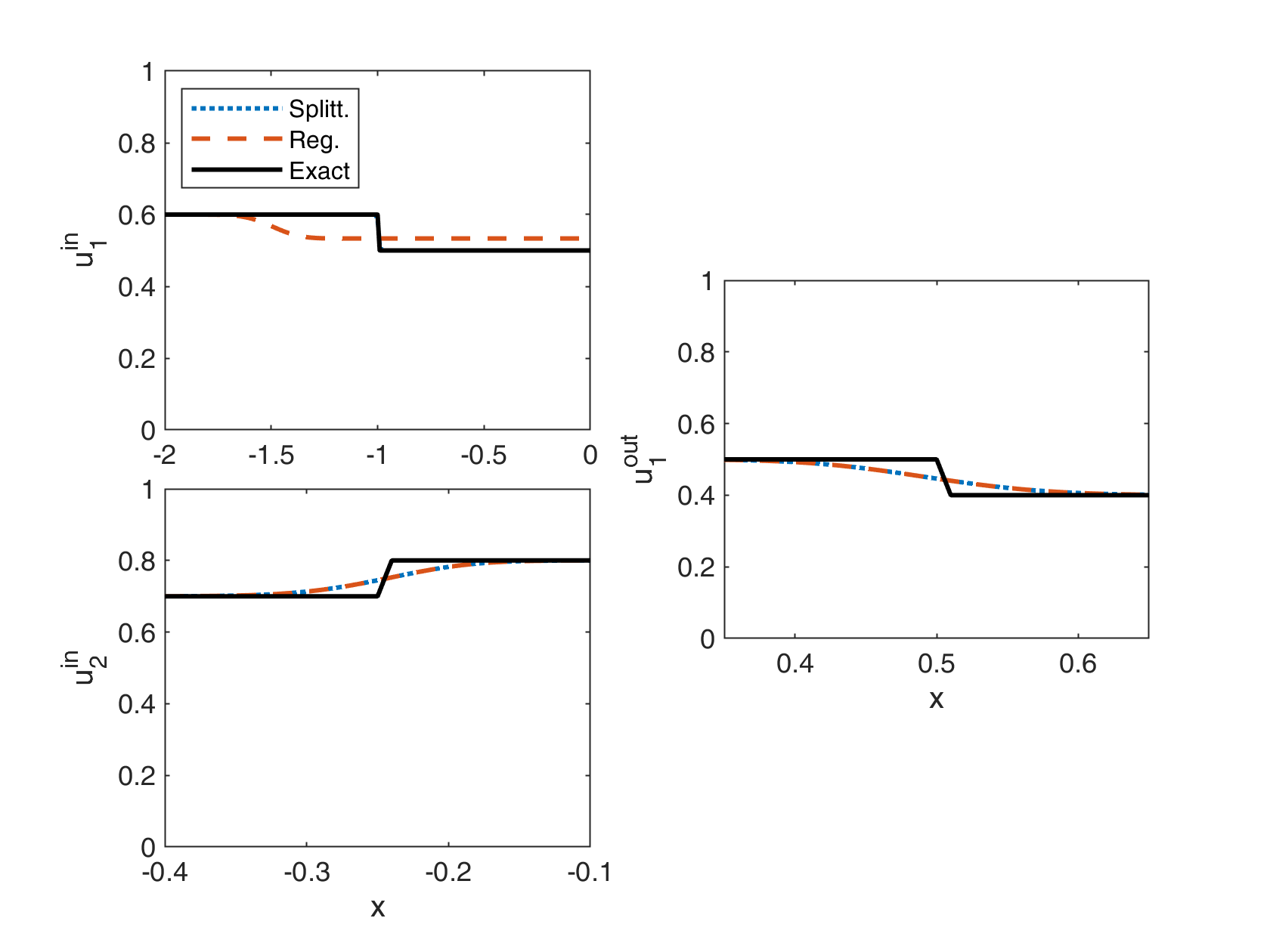}
   \caption{Second example with $\lambda=0.1$ and $\Delta x=0.01$.}   
    \label{fig:2to1}
\end{figure}

For  $\lambda=0.1$ and $\Delta x=0.01$ the approximate solutions and the exact solution of the second example are displayed in figure \ref{fig:2to1}.
In particular, on the first incoming road the splitting algorithm approximates the solution much better.
\rv{For the considered resolution the difference between the exact solution and its approximation by the splitting algorithm is not visible on the first road.}
On the other roads the approximate solutions are very similar. Here, we zoomed in to display the differences to the exact solution.

\section{Conclusion}
We have presented a Riemann solver at a junction for conservation laws with discontinuous flux.
We have adapted the splitting algorithm of \cite{Towers} to networks and demonstrated its validity in comparison with the exact solution. We have also pointed out that the splitting algorithm on networks is faster and more accurate than the approach with a regularized flux.
Future research includes the investigation of other network models, where the flux is discontinuous.

\section*{Acknowledgment}
J.F. was supported by the German Research Foundation (DFG) under grant HE 5386/18-1 and
S.G.~under grant GO 1920/10-1.

\bibliographystyle{siam}
\bibliography{Literatur}

\end{document}